\documentclass[draft]{amsart}
\usepackage{amssymb,amsxtra,verbatim}

\theoremstyle{plain}
\newtheorem{thm}[equation]{Theorem} 
\newtheorem{lem}[equation]{Lemma} 
\newtheorem{cor}[equation]{Corollary}
\newtheorem{prop}[equation]{Proposition} 
\newtheorem{example}[equation]{Example}
\newtheorem{examples}[equation]{Examples}

\theoremstyle{definition}
\newtheorem{defn}[equation]{Definition}
\newtheorem{notation}[equation]{Notation}

\newtheorem{rem}[equation]{Remark}

\numberwithin{equation}{section}

\newcommand{\Q}{{\mathbb{Q}}}
\newcommand{\N}{{\mathbb{N}}}
\newcommand{\C}{{\mathbb{C}}}
\newcommand{\R}{{\mathbb{R}}}

\newcommand{\Z}{{\mathbb{Z}}}

\newcommand{\K}{{\mathbb{K}}}

\newcommand{\calC}{{\mathcal{C}}}
\newcommand{\cO}{{\mathcal{O}}}

\newcommand{\calB}{{\mathcal{B}}}
\newcommand{\calD}{{\mathcal{D}}}\newcommand{\cD}{{\mathcal{D}}}

\newcommand{\calM}{{\mathcal{M}}}\newcommand{\cM}{{\mathcal{M}}}
\newcommand{\calN}{{\mathcal{N}}}

\newcommand{\calS}{{\mathcal{S}}}

\DeclareMathAlphabet{\euls}{U}{eus}{m}{n}

\DeclareMathAlphabet{\eulr}{U}{eur}{m}{n}
\DeclareMathAlphabet{\eulrb}{U}{eur}{b}{n}

\newcommand{\eF}{{\euls{F}}}

\newcommand{\eI}{{\euls{I}}}

\newcommand{\eL}{{\euls{L}}}

\newcommand{\fa}{{\mathfrak{a}}}
\newcommand{\g}{{\mathfrak{g}}}

\newcommand{\fh}{{\mathfrak{h}}}
\newcommand{\h}{{\mathfrak{h}}}

\newcommand{\fn}{{\mathfrak{n}}}
\newcommand{\fp}{{\mathfrak{p}}}

\newcommand{\fr}{{\mathfrak{r}}}

\newcommand{\fsl}{{\mathfrak{sl}}}

\newcommand{\ft}{{\mathfrak{t}}}
\newcommand{\fu}{{\mathfrak{u}}}
\newcommand{\fs}{{\mathfrak{s}}}

\newcommand{\mbx}{{\mathbf{x}}}
\newcommand{\mby}{{\mathbf{y}}}

\newcommand{\mbB}{{\mathbf{B}}}

\newcommand{\mbO}{{\mathbf{O}}}

\newcommand{\mbV}{{\mathbf{V}}}

\newcommand{\mbX}{{\mathbf{X}}}
\newcommand{\mbY}{{\mathbf{Y}}}
\newcommand{\mbZ}{{\mathbf{Z}}}

\def\preisomto{\vbox{\hbox to
                 15pt{\hfill$\sim$\hfill}\nointerlineskip\vskip
                 -0.3pt 
                 \hbox to 16pt{\rightarrowfill}}}

\def\isomto{\mathop{\preisomto}}

\def\prelongisomto{\vbox{\hbox to
                18pt{\hfill$\sim$\hfill}\nointerlineskip\vskip -0.3pt
                \hbox to 19pt{\rightarrowfill}}}  

\def\longisomto{\mathop{\prelongisomto}}

\newcommand{\doublelongrightarrow}{\longrightarrow \kern-14pt
\longrightarrow}

\newcommand{\sto}{\twoheadrightarrow}
\newcommand{\ito}{\hookrightarrow}
\newcommand{\lto}{\longrightarrow}

\def\trait{\hbox to 4mm{\hrulefill}}

\newcommand{\limply}{\Longrightarrow}

\def\Aut{\operatorname {Aut}}

\def\Ker{\operatorname {Ker}}
\def\Im{\operatorname {Im}}

\def\SL{\operatorname {SL}}
\def\Ext{\operatorname {Ext}}

\def\grade{\operatorname {j}}

\def\gr{\operatorname {gr}}
\def\rk{\operatorname {rk}}
\def\codim{\operatorname {codim}}

\def\GK{\operatorname {GKdim}}

\def\Lie{\operatorname {Lie}}
\def\Der{\operatorname {Der}}

\def\ord{\operatorname {ord}}
\def\trdeg{\operatorname {trdeg}}
\def\coh{\operatorname {H}}

\newcommand{\dual}[2]{{\langle {#1} , {#2} \rangle}}

\newcommand{\vpi}{\varpi}
\newcommand{\vphi}{\varphi}

\newcommand{\ges}{\geqslant}

\newcommand{\Fo}{{F_{w_0}}}
\newcommand{\Deltacheck}{{\Delta\spcheck}}
\newcommand{\Deltapc}{{\Delta_{+}\spcheck}}
\newcommand{\alphacheck}{{\alpha\spcheck}}

\newcommand{\alphac}{{\alpha\spcheck}}
\newcommand{\gammas}{{\gamma^*}}
\newcommand{\Ygamma}{\mbY_\gamma}
\newcommand{\YGamma}{\mbY_\Gamma}
\newcommand{\ombY}{{\overline{\mbY}}}

\newcommand{\DYg}{{\cD(\mbY_\gamma)}}
\newcommand{\OX}{{\cO(\mbX)}}
\newcommand{\OXw}{{\cO(\mbX)^w}}
\newcommand{\DX}{{\cD(\mbX)}}
\newcommand{\sOX}{{\cO_\mbX}}
\renewcommand{\L}[1]{{\eL({#1})}}
\newcommand{\twoslash}{{{/} \hskip -3pt {/}}}
\newcommand{\wh}{{\widehat{\fh}}}

\begin{document}
   

  
\title[Differential Operators on the Basic Affine
Space]{Differential Operators and Cohomology Groups on the
  Basic Affine Space} \author{T.~Levasseur}
\address{D\'epartement de Math\'ematiques, Universit\'e de
  Brest, 29238 Brest cedex~3, France}
\email{Thierry.Levasseur@univ-brest.fr} \author{J.~T.~Stafford}
\address{Department of Mathematics, University of Michigan, Ann
  Arbor, MI 48109-1109, USA} \email{jts@umich.edu} \thanks{The
  second author was supported in part by the NSF through the
  grants DMS-9801148 and DMS-0245320. Part of this work was
  done while he was visiting the Mittag-Leffler Institute and
  he would like to thank the Institute for its financial
  support and hospitality.}   \dedicatory{This paper is
  dedicated to Tony Joseph on the occasion of his
  $60^{\text{th}}$ birthday.}
   \begin{abstract}  
     We study the ring of differential operators $\cD(\mbX)$ on
     the basic affine space $\mbX=G/U$ of a complex semisimple
     group $G$ with maximal unipotent subgroup $U$.  One of the
     main results shows that the cohomology group
     $\mathrm{H}^*(\mbX,\sOX)$ decomposes as a finite direct
     sum of non-isomorphic simple $\cD(\mbX)$-modules, each of
     which is isomorphic to a twist of $\cO(\mbX)$ by an
     automorphism of $\cD(\mbX)$.
    
     We also use $\cD(\mbX)$ to study the properties of
     $\cD(\mbZ)$ for highest weight varieties $\mbZ$. For
     example, we prove that $\mbZ$ is $\cD$-simple in the sense
     that $\cO(\mbZ)$ is a simple $\cD(\mbZ)$-module and
     produce an irreducible $G$-module of differential
     operators on $\mbZ$ of degree $-1$ and specified order.
   \end{abstract}
   
   \maketitle

\tableofcontents

\clearpage


\section{Introduction}\label{intro}

Fix a complex semisimple, connected and simply connected Lie
group $G$ with maximal unipotent subgroup $U$ and Lie algebra
$\g$.  Then the \emph{basic affine space} is the quasi-affine
variety $\mbX=G/U$.  The ring of global differential operators
$\cD(\mbX)$ has a long history, going back to the work
\cite{GK} of Gelfand and Kirillov in the late 60's who used
this space to formulate and partially solve their conjecture
that the quotient division ring of the enveloping algebra
$U(\g)$ should be isomorphic to a Weyl skew field.

The variety $\mbX$ is only quasi-affine and, when $\g$ is not
isomorphic to a direct sum of copies of $\mathfrak{sl}(2)$, the
affine closure $\overline{\mbX}$ of $\mbX$ is singular.  In
general, rings of differential operators on a singular variety
$\mbZ$ can be quite unpleasant; for example, and in contrast to
the case of a smooth affine variety, $\cD(\mbZ)$ need not be
noetherian, finitely generated or simple and (conjecturally) it
will not be generated by the derivations $\mathrm{Der}_{\mathbb
  C}(\mbZ)$.  Moreover, the canonical module $\cO(\mbZ)$ need
not be a simple $\cD(\mbZ)$-module.  Recently, Bezrukavnikov,
Braverman and Positselskii \cite{BBP} proved a remarkable
result on the structure of $\cD(\mbX)$ which shows that it
actually has very pleasant properties.  Specifically, they
proved that there exist automorphisms $\{F_w\}_w$, indexed by
the Weyl group $W$ of $G$, such that for any nonzero
$\cD(\mbX)$-module $M$ there exists $w\in W$ such that
$\cD_\mbX\otimes_{\cD(\mbX)}M^{w}\not=0$, where $M^w=M^{F_w}$
is the twist of $M$ by $F_w$.  (The $F_w$ should be thought of
as analogues of partial Fourier transforms.)  Since $\mbX$ is
smooth this implies that, for any finite open affine cover
$\{\mbX_i\}_i$ of $\mbX$, the ring
$\bigoplus_{i,w}\cD(\mbX_i)^{w}$ is a noetherian, faithfully
flat overring of $\cD(\mbX)$.  As is shown in \cite{BBP}, it
follows easily that $\cD(\mbX)$ is a noetherian domain of
finite injective dimension.

The aim of this paper is to extend and apply the results of
\cite{BBP}.  Our first result, which combines
Proposition~\ref{cor26}, Theorem~\ref{thm29} and
Theorem~\ref{thm36}, further elucidates the structure of
$\cD(\mbX)$.

\begin{prop}\label{basicprops-intro}
  Let $\mbX=G/U$ denote the basic affine space of $G$.  Then:
  \begin{enumerate}
  \item $\cD(\mbX)$ is a simple ring satisfying the
    Auslander-Gorenstein and Cohen-Mac\-aulay conditions (see
    Section~\ref{structure-section} for the definitions);
  \item $\cD(\mbX)$ is (finitely) generated, as a $\mathbb
    C$-algebra, by $\{\cO(\mbX)^w : w\in W\}\cup
    \widehat{\fh}$.
 \end{enumerate}
 \end{prop}
 
 Note that $\cD(\mbX)$ is quite a subtle ring; for example, if
 $G=\SL(3,\mathbb C)$ then $\cD(\mbX) \cong
 U(\mathfrak{so}(8))/J$, where $J$ is the {\em Joseph ideal} as
 defined in~\cite{Jos1} (see Example~\ref{sl-3}).

 The variety $\mbX=G/U$ has a natural left action of $G$ and a
 right action of the maximal torus $H$ for which $B=HU$ is a
 Borel subgroup.  Differentiating these actions gives
 embeddings of $\g=\mathrm{Lie}(G)$, respectively
 $\widehat{\fh}=\mathrm{Lie}(H)$ into $\mathrm{Der}_{\mathbb
   C}(\mbX)$. It follows easily from the simplicity of
 $\cD(\mbX)$ that $\cO(\mbX)=\mathrm{H}^0(\mbX,\sOX)$ is a
 simple $\cD(\mbX)$-module. One of the main results of this
 paper extends this to describe the $\cD(\mbX)$-module
 structure of the full cohomology group
 $\mathrm{H}^*(\mbX,\sOX)$:

\begin{thm}\label{cohom-intro}\emph{[Theorem~\ref{thm57}]}
  Let $\mbX=G/U$ and set $\cM=\cD(\mbX)/\cD(\mbX)\g$.
  \begin{enumerate}
  \item $\cM\cong \mathrm{H}^*(\mbX,\sOX)$.
  \item For each $i$, $\mathrm{H}^i(\mbX,\sOX)\cong \bigoplus
    \Bigl\{\cO(\mbX)^w : w\in W : \ell(w)=i\Bigr\}$.
  \item For $w\not=v\in W$, the $\cD(\mbX)$-modules
    $\cO(\mbX)^w$ and $\cO(\mbX)^{v}$ are simple and
    nonisomorphic.
\end{enumerate}
\end{thm}

A result analogous to Theorem~\ref{cohom-intro}, but in the
$l$-adic setting, has been proved in \cite[Lemma~12.0.1]{Po}.
It is not clear to us what is the relationship between the two
results.

A key point in the proof of Theorem~\ref{cohom-intro} is that,
by the Borel-Weil-Bott theorem, one has an explicit description
of the $G$-module structure of the $ \mathrm{H}^*(\mbX,\sOX)$
and one then proves Theorem~\ref{cohom-intro} by comparing that
structure with the $G$-module structure of the $\cO(\mbX)^w$.
Theorem~\ref{cohom-intro} is rather satisfying since it relates
the left ideal of $\cD(\mbX)$ generated by the derivations
coming from $\g$ to the only ``obvious'' simple
$\cD(\mbX)$-modules $\cO(\mbX)^w$.  In contrast, if one
considers the left ideal generated by all the ``obvious''
derivations, $\cD(\mbX)\g+\cD(\mbX)\wh$, then one obtains:

\begin{prop}\label{deriv-intro}\emph{[Theorem~\ref{thm54}]}  
  As left $\cD(\mbX)$-modules,
  $$
  \frac{\cD(\mbX)}{\cD(\mbX)\g +\strut\cD(\mbX)\wh} \ = \ 
  \frac{\cD(\mbX)}{\cD(\mbX)\mathrm{Der}_\C(\mbX)} \ \cong \ 
  \cO(\mbX).
  $$
\end{prop}
 
This proposition is somewhat surprising since, at least when
$\g$ is not a direct sum of copies of $\fsl(2)$, one can show
that $\cD(\mbX)$ is \emph{not} generated by $\cO(\mbX)$ and
$\mathrm{Der}_\C(\mbX)$ as a $\C$-algebra (see
Corollary~\ref{cor61}). Of course, the analogue of
Proposition~\ref{deriv-intro} for smooth varieties is standard.

A natural class of varieties associated to $G$ are
\emph{S-varieties}: closures $\overline{\mbY}$ of a $G$-orbit
$\mbY=\mbY_\Gamma=G.v_\Gamma$ where $v_\Gamma$ is a sum of
highest weight vectors in some finite dimensional $G$-module.
In such a case there exists a natural surjection $
\mbX\twoheadrightarrow \mbY=G/S_\Gamma$, for the isotropy group
$S_\Gamma$ of $v_\Gamma$. This induces, by restriction of
operators, a map $\psi_\Gamma : \cD(\mbX)^{S_\Gamma}\to
\cD(\overline{\mbY})$ and allows us to use our structure
results on $\cD(\mbX)$ give information on $\cD(\mbY)$. For
this to be effective we need the mild assumption that
$\overline{\mbY}$ is normal and
$\mathrm{codim}_{\overline{\mbY}}(\overline{\mbY}\smallsetminus
\mbY) \geq 2$ or, equivalently, that $\cO(\mbY)=\cO(\ombY)$
(see Theorem~\ref{vinberg-popov} for further equivalent
conditions).

\begin{cor}
\label{wt-var-intro}\emph{[Proposition~\ref{thm33'}]}
Let $\ombY$ be an S-variety such that $\cO(\mbY)=\cO(\ombY)$.
Then $\cO(\mbY) $ is a simple $\cD(\mbY)$-module.
\end{cor}

When $\ombY=\ombY_\Gamma$ is singular, this result says that
there exists operators $D\in \cD(\mbY)$ that cannot be
constructed from derivations; these are ``exotic'' operators in
the terminology of \cite{AB}).  As the name suggests, exotic
operators can be hard to construct---see \cite{AB} or
\cite{BK}, for example---but in our context their construction
is easy; they arise as $\psi_\Gamma
F_{w_0}(\cO(\mbY_{\Gamma^*}))$.  In the special case of a
highest weight variety (this is just an S-variety
$\ombY_\gamma=\ombY_\Gamma$ for $\Gamma=\mathbb{N}\gamma$) we
can be more precise about these operators.

\begin{cor}\label{exotic-intro}  \emph{[Corollary~\ref{exotic-cor}]}
  Suppose that $\ombY=\ombY_\gamma$ is a highest weight
  variety.  Then there exists an irreducible $G$-module $E\cong
  V(\gamma)$ of differential operators on $\mbY$ of degree $-1$
  and order $\dual{\gamma}{2\rho\spcheck}$.
 \end{cor}
 
 Corollary~\ref{wt-var-intro} also proves the Nakai conjecture
 for S-varieties satisfying the hypotheses of that result and
 this covers most of the known cases of normal, singular
 varieties for which the conjecture is known.  See
 Section~\ref{sec3} for the details.

 A fundamental question in the theory of differential operators
 on invariant rings asks the following. Suppose that $Q$ is a
 (reductive) Lie group acting on a finite dimensional vector
 space $V$ such that the fixed ring $\cO(V)^Q$ is singular.
 Then, is the natural map $\cD(V)^Q\to \cD(\cO(V)^Q)$
 surjective?  A positive answer to this question is known in a
 number of cases and this has had significant applications to
 representation theory; see, for example,
 \cite{Jos2,LS0,LS1,Sc} and Remark~\ref{lss-comment}.  It is
 therefore natural to ask when the analogous map $\psi_\Gamma :
 \cD(\mbX)^{S_\Gamma}\to \cD(\mbY_\Gamma)$ is surjective.
 Although we do not have a general answer to this question, we
 suspect that $\psi_\Gamma$ will usually not be surjective. As
 evidence for this we prove that this is the case for one of
 the fundamental examples of an S-variety: the closure of the
 minimal orbit in a simple Lie algebra~$\g$.
 
\begin{prop} \label{nonsurj-intro}  
\label{intro-minorbit}\emph{[Theorem~\ref{non-exotic}]}
If $\mbY_\gamma $ is the minimal (nonzero) nilpotent orbit of a
simple classical Lie algebra $\g$, then $\psi_\gamma$ is
surjective if and only if $\g=\mathfrak{sl}(2)$ or
$\g=\mathfrak{sl}(3)$.
\end{prop}
  
  \noindent
  {\bf Acknowledgement.} We would like to thank the referee for
  some helpful comments.
 

\section{Preliminaries} 
\label{sec1}

In this section we describe the basic results and notation we
need from the literature, notably the relevant results from
\cite{BBP}.  The reader is also referred to \cite{GK, HV, Sa0,
  Sa1} for the interrelationship between differential operators
on the base affine space and the corresponding enveloping
algebra.

We begin with some necessary notation.  The base field will
always be the field $\C$ of complex numbers.  Let $G$ be a
connected simply-connected semisimple algebraic group of rank
$\ell$, $B$ a Borel subgroup, $H \subset B$ a maximal torus, $U
\subset B$ a maximal unipotent subgroup of $G$.  The Lie
algebra of an algebraic group is denoted by the corresponding
gothic character; thus $\g = \mathrm{Lie}(G)$, $\fh =
\mathrm{Lie}(H)$ and $\fu = \mathrm{Lie}(U)$.

We will use standard Lie theoretic notation, as for example
given in \cite{Bou}.  In particular, let $\Delta$ denote the
root system of $(\g,\fh)$ and fix a set of positive roots
$\Delta_+$ such that $\fu = \bigoplus_{\alpha \in \Delta_+}
\fu_\alpha$. Denote by $W$ the Weyl group of $\Delta$.  Let
$\Lambda$ be the weight lattice of $\Delta$ which we identify
with the character group of $H$. The set of dominant weights is
denoted by $\Lambda^+$. Fix a basis $\Sigma =
\{\alpha_1,\dots,\alpha_\ell\}$ of $\Delta_+$ and write
$\{\vpi_1,\dots,\vpi_\ell\}$ for the fundamental weights; thus
$\dual{\vpi_i}{\alpha_j\spcheck} = \delta_{ij}$. Denote by
$\rho$, respectively $\rho\spcheck$, the half sum of the
positive roots, respectively coroots.  Let $w_0$ be the longest
element of $W$ and set $\lambda^* = -w_0(\lambda)$ for all
$\lambda \in \Lambda$.  Then $ w_0(\rho) = -\rho$,
$w_0(\rho\spcheck)=-\rho\spcheck$ and, by \cite[Chapter~6,
1.10, Corollaire]{Bou},
$$
\dual{\lambda^*}{\rho\spcheck}=\dual{\lambda}{\rho\spcheck}
= \sum_{i=1}^\ell m_i \qquad\text{if}\quad \lambda =
\sum_{i=1}^\ell m_i \alpha_i.
$$

For each $\omega \in \Lambda^+$ let $V(\omega)$ be the
irreducible $G$-module of highest weight $\omega$ and
$V(\omega)_\mu$ be the subspace of elements of weight $\mu \in
\Lambda$; we will denote by $v_\omega$ a highest weight vector
of $V(\omega)$. Recall that the $G$-module $V(\omega)^*$
identifies naturally with $V(\omega^*)$.  If $E$ is a locally
finite $(G \times H)$-module, we denote by $E[\lambda]$ the
isotypic $G$-component of type $\lambda \in \Lambda^+$ and by
$E^\mu$ the $\mu$-weight space for the action of $H \equiv
\{1\} \times H$. Hence,
$$
E = \bigoplus\left\{ E[\lambda]^\mu : \mu \in \Lambda,\,
  \lambda \in \Lambda^+\right\}.
$$

Let $\mbY$ be an algebraic variety. We denote by $\cO_{\mbY}$
the structural sheaf of $\mbY$ and by $\cO(\mbY)$ its algebra
of regular functions. The sheaf of differential operators on
$\mbY$ is denoted by $\cD_\mbY$ with global sections
$\cD(\mbY)=\text{H}^0(\mbY,\cD_\mbY )$.  The $\cO(\mbY)$-module
of elements of order $\le k$ in $\cD(\mbY)$ is denoted by
$\cD_k(\mbY)$, and the order of $D \in \cD(\mbY)$ will be
written $\ord D$. Now suppose that $\mbY$ is an irreducible
quasi-affine algebraic variety, embedded as an open subvariety
of an affine variety $\overline{\mbY}$.  We will frequently use
the fact that, if $\overline{\mbY}$ is normal with
$\codim_{\overline{\mbY}}(\overline{\mbY} \smallsetminus \mbY)
\ge 2$, then $\cO(\mbY)= \cO(\overline{\mbY})$ and so
$\cD(\mbY)= \cD(\overline{\mbY})$.

Let $Q$ be an affine algebraic group. We say that an algebraic
variety $\mbY$ is a $Q$-variety if it is equipped with a
rational action of $Q$.  For such a variety, $\cO(\mbY),
\cD_k(\mbY)$ and $\cD(\mbY)$ are locally finite $Q$-modules,
with the action of $a \in Q$ on $\vphi \in \cO(\mbY)$ and $D
\in \calD(\mbY)$ being defined by $a.\vphi(y) =
\vphi(a^{-1}.y)$ for $y \in \mbY$, respectively $(a.D)(\vphi) =
a.D(a^{-1}.\vphi)$.  If $\mbY$ is also an affine variety such
that $\cO(\mbY)^Q$ is an affine algebra, we define the
\emph{categorical quotient} $\mbY \twoslash Q$ by $\cO(\mbY
\twoslash Q)= \cO(\mbY)^Q$. We then have a restriction morphism
$$
\psi : \cD(\mbY)^Q \lto \cD(\mbY \twoslash Q), \qquad
\psi(D)(f) = D(f) \ \, \text{for $f \in \cO(\mbY)^Q$.}
$$
Notice that $\psi(\cD_k(\mbY)^Q) \subseteq \cD_k(\mbY
\twoslash Q)$ for all $k$. In many cases $\cO(\mbY)$ will be a
$\mathbb Z$-graded algebra, $\cO(\mbY)=\bigoplus_{n\in \mathbb
  Z} \cO^n$, in which case $\cD(\mbY)$ has an induced $\mathbb
Z$-graded structure, with the $n^{\text{th}}$ graded piece
being
\begin{equation}\label{D-graded}
\cD(\mbY)^n=\{\theta\in \cD(\mbY) : \theta(\cO^r)\subseteq \cO^{r+n}
\ \text{for all}\ r\in \mathbb Z\}.
\end{equation}
In this situation, $\cD(\mbY)^n$ will be called operators of
\emph{degree} $n$.

Assume now that $V$ is a $(G \times H)$-module and that $\mbY$
is a $(G \times H)$-subvariety of~$V$. Then
$$
\cO(\mbY) = \bigoplus_{\substack{{\scriptstyle\mu \in \Lambda}\\
    {\scriptstyle\lambda \in \Lambda^+ }} }
\cO(\mbY)[\lambda]^\mu \ \subset \ 
\cD(\mbY) = \bigoplus_{\substack{{\scriptstyle\mu \in \Lambda}\\
    {\scriptstyle\lambda \in \Lambda^+ }} }
\cD(\mbY)[\lambda]^\mu.
$$
It follows easily from the definitions that
\begin{equation}\label{D-mu}
\cD(\mbY)^\mu = \bigl\{d \in \cD(\mbY) :
d(\cO(\mbY)^\lambda) \subseteq \cO(\mbY)^{\lambda + \mu} \ 
\text{for all $\lambda \in \Lambda$}\bigr\}.
\end{equation}
One clearly has a surjection $S(V^*)[\lambda]^\mu \sto
\cO(\mbY)[\lambda]^\mu$. Furthermore, if $V=V^{-\gammas}$ for
some $0 \ne \gammas \in \Lambda$, we will identify $S^m(V^*)$
with $S(V^*)^{m\gammas}$ and obtain the surjective $G$-morphism
$S^m(V^*)[\lambda] \sto \cO(\mbY)[\lambda]^{m\gammas}$.  In
this case, $\cO(\mbY) = \bigoplus_{m\in \N}
\cO(\mbY)^{m\gammas}$ and $\cD(\mbY) = \bigoplus_{m\in \Z}
\cD(\mbY)^{m\gammas}$ are $\mathbb Z$-graded and \eqref{D-mu}
can be interpreted as saying that $\cD(\mbY)^{m\gammas}$ is the
set of differential operators of degree $m$ on $\mbY$.

The previous results apply in particular to the basic affine
space $\mbX = G/U$.  We need to collect here a few facts about
the $(G \times H)$-variety $\mbX$ and its canonical affine
closure $\overline{\mbX}$.  For these assertions, see, for
example, \cite{GK, Gr2, VP}.  Set $V = V(\vpi_1) \bigoplus
\cdots \bigoplus V(\vpi_\ell)$ and recall that there is an
isomorphism
$$
\mbX \longisomto G.(v_{\vpi_1} \oplus \cdots \oplus
v_{\vpi_\ell}) \subset V, \qquad \bar{g} \mapsto g.(v_{\vpi_1}
\oplus \cdots \oplus v_{\vpi_\ell}),
$$
where $\bar{g}$ denotes the class of $g \in G$ modulo $U$.
We identify $\mbX$ with this $G$-orbit and write
$$
\cO = \cO(\mbX) \quad\text{ and }\quad \cD = \cD(\mbX).
$$
Then the Zariski closure $\overline{\mbX}$ of $\mbX$ in $V$
is a normal irreducible affine variety that satisfies
$\codim_{\overline{\mbX}}(\overline{\mbX} \smallsetminus \mbX)
\ge 2$; thus $\cO=\cO(\overline{\mbX})$, etc.

Identify $w_0$ with an automorphism of $H$; thus $\lambda^*(h)
= \lambda(w_0(h^{-1}))$ for all $\lambda \in \Lambda$ and $ h
\in H$.  Define the twisted (right) action of $H$ on $\mbX$ by
$r_h.\bar{g} = \overline{gw_0(h)}$. Endow the $G$-module $V$
with the action of $H \equiv \{1\} \times H$ defined by
$$r_h.(u_1 \oplus \cdots \oplus u_\ell) = \vpi_1^*(h^{-1}) u_1
\oplus \cdots \oplus \vpi_\ell^*(h^{-1}) u_\ell.$$
Then the
embedding $\mbX \ito V$ is a morphism of $(G \times
H)$-varieties.  The induced (left) action of $H$ on $\cO$ is
then given by $(r_h.f)(\bar{g})=f(r_{h^{-1}}.\bar{g})$ for all
$f \in \cO, h \in H, \bar{g} \in \mbX$.  It follows that
$\cO[\nu^*] = \cO^{\nu^*}$, $\cO^\mu= 0$ if $\mu \notin
\Lambda^+$, and we can decompose the $(G\times H)$-module $\cO$
as:
\begin{equation}\label{O-decomp}
 \cO = \bigoplus_{\lambda \in \Lambda^+}
  \cO^\lambda, \qquad \cO^\lambda = \cO[\lambda] \cong
  V(\lambda).
\end{equation}
The final isomorphism in \eqref{O-decomp} comes from the fact
that the $G$-action on $\cO$ is multiplicity free
\cite[Theorem~2]{VP}.  Notice that the algebra $\cO$ is
(finitely) generated by the $G$-modules $\cO^{\vpi_j}$, $1 \le
j \le \ell$.  Also, \eqref{D-mu} implies that
$\cD^\mu(\cO^\lambda)= 0$ when $\lambda + \mu$ is not dominant.

\begin{notation}\label{h-hat}
  The differentials of the actions of $G$ and $H$ (via $h
  \mapsto r_h$) on $\mbX$ yield morphisms of algebras $\imath :
  U(\g) \to \cD$ and $\jmath : U(\fh) \to \cD$. By
  \cite[Corollary~8.1 and Lemma~9.1]{GK}, $\imath$ and $\jmath$
  are injective; from now on we will identify $U(\g)$ with
  $\imath(U(\g))$ but write $\widehat{\h}=\jmath(\h)$ and
  $U(\widehat{\h})=\jmath(U(\h))$, to distinguish these objects
  from their images under $\imath$.
\end{notation}

Let $0 \ne M$ be a left $\cD$-module and $\tau \in \Aut(\cD)$,
the $\mathbb C$-algebra automorphism group of $\cD$. Then the
\emph{twist of $M$ by $\tau$} is the $\cD$-module $M^\tau$
defined as follows: $M^\tau = M$ as an abelian group but $a
\cdot x = \tau(a)x$ for all $a \in \cD, x \in M$.  Recall also
that the localization of $M$ on $\mbX$ is
$$
L(M) = \cD_{\mbX} \otimes_{\cD} M \cong \cO_\mbX \otimes_\cO
M.
$$
Thus $L(M)$ is a quasi-coherent left $\cD_\mbX$-module.
Clearly $L(M)=0$ when $M$ is supported on $\overline{\mbX}
\smallsetminus \mbX$ but, remarkably, one can obtain a nonzero
localization by first twisting the module $M$:

\begin{thm} \cite[Proposition~3.1 and Theorem~3.4]{BBP}
\label{thm21}        Let $\cD=\cD(\mbX)$. Then:
\begin{enumerate}
\item[(1)] There exists an injection of groups $F : W \ito
  \Aut(\cD)$ written $w \mapsto F_w$.
\item[(2)] For each $\cD$-module $M \ne 0$ there exists $w\in
  W$ such that $L(M^{F_w})~\ne~0$.  \qed
\end{enumerate}
\end{thm}

This theorem is also valid for right modules.  The morphisms
$F_w$ can be regarded as variants of partial Fourier
transforms, and more details on their structure and properties
can be found in \cite{BBP}.

It will be convenient to reformulate Theorem~\ref{thm21}, for
which we need some notation.  Since $\mbX$ is an open subset of
the non singular locus of $\overline{\mbX}$, we can fix an open
(and smooth) affine cover of $\mbX$:
\begin{equation}\label{cover-equ}
\mbX = \bigcup_{j=1}^k \mbX_j, \quad \text{where} \quad \mbX_j
= \{x \in 
\overline{\mbX} : f_j(x) \ne 0\},
\end{equation}
for the appropriate $f_i\in \cO(\mbX)$.  Notice that
$\cD(\mbX_j) = \cD[f_j^{-1}]$ for each $j$ and so $\cD(\mbX_j)$
is a flat $\cD$-module.  If $M$ is $\cD$-module and $w\in W$,
let $M^{w}$ denote the twist of $M$ by $F_{w}$.

For each pair $(w,j)$ with $w\in W$ and $1 \le j \le k$, we
have an injective morphism of algebras $\phi_{wj} : \cD \ito
\cD(\mbX_j)$ given by $\phi_{wj}(d) = F_{w}^{-1}(d)$. We write
$\cD(\mbX_j)_{w}$ for the algebra $\cD(\mbX_j)$ regarded as an
overring of $\cD$ under this embedding.  The significance of
this construction is that, for any left $\cD$-module $M$, the
map $ d\otimes v\mapsto d\otimes v$ induces an isomorphism $
\cD(\mbX_j) \otimes_\cD M^{{w}} \cong \cD(\mbX_j)_{w}
\otimes_\cD M$ of left $\cD(\mbX_j) $-modules.
Theorem~\ref{thm21} can now be rewritten as follows.

\begin{cor}
\label{cor24}
Let $0 \ne M$ be a left $\cD$-module. Then:
\begin{enumerate}
\item There exists a pair $(w,j)$ such that $\cD(\mbX_j)_{w}
  \otimes_\cD M \ne 0$.
\item Set $R_w =\bigoplus_{j=1}^k \cD(\mbX_j)_{w}$ and
  $R=\bigoplus_{w\in W} R_w$.  Then $R$ is a faithfully flat
  (left or right) overring of $\cD$.
\end{enumerate}\end{cor}

\begin{proof}
  Part (1) is just a reformulation of Theorem~\ref{thm21}. This
  in turn implies that $R$ is a faithful right $\cD$-module. It
  is flat since $\cD(\mbX_j)_{w}$ is isomorphic, as a
  $\cD$-module, to the localization of $\cD$ at the powers of
  $F_w(f_j)$.  Since Theorem~\ref{thm21} also holds for right
  modules, the same argument shows that $R$ is a faithful flat
  left $\cD$-module.
\end{proof}

When $\g=\mathfrak{sl}(2)$ it is easy to check that
$\overline{\mbX}=\mathbb C^2$, and so there is no subtlety to
the structure of either $ \overline{\mbX}$ or $\cD$. However,
when $\g\not=\mathfrak{sl}(2)^m$, $\overline{\mbX}$ will be
singular and $\cD$ will be rather subtle. The simplest example
is:

\begin{example}
\label{sl-3}
Assume that $\g=\mathfrak{sl}(3)$ and set $\mbX=\SL(3)/U$.
Then $\overline{\mbX}$ is the quadric $\sum_{i=1}^3 u_iy_i=0$
inside $\C^6$. Moreover, $\cD=\cD(\mbX) \cong
U(\mathfrak{so}(8))/J$, where $J$ is the Joseph ideal.

An explicit set of generators of $\cD$ is given in
\cite[(2.2)]{LS2}.  The algebra $F_{w_0}(\cO) $ is generated by
the operators $\{\Phi_j,\Theta_j\}$ of order $2$ from
\cite[(2.2)]{LS2}.
\end{example}

\begin{proof}
  The proof of the first assertion is an elementary and
  classical computation, which is left to the reader.  The
  second assertion then follows from \cite[Remark~3.2(v) and
  Corollary~A]{LSS}.  The claim in the second paragraph will
  not be used in this paper and so is left to the interested
  reader (the next proposition may prove useful). A second way
  of interpreting this example is given in
  Remark~\ref{lss-comment}.
\end{proof}

In the main body of the paper we will need some more technical
results from \cite{BBP} about the automorphisms $F_w$ and for
the reader's convenience we record them in the following
proposition. This summarizes Lemma~3.3, Corollary~3.10,
Proposition~3.11 and the proof of Lemma~3.12 of that paper.

\begin{prop}
\label{prop22}
Let $\eta \in \Lambda^+$, $\mu \in \Lambda$ and $w \in W$.
Then:
\begin{enumerate}
\item[(1)] $F_w$ is $G$-linear and $F_w(\cD^\mu) =
  \cD^{w(\mu)}$. Thus $F_w(\cO^\eta) \cong V(\eta)$ as
  $G$-modules.
\item[(2)] $F_w(h) = w(h) + \dual{w(h) -h}{\rho}$ for all $h
  \in \wh$.
\item[(3)] $\ord F_w(d) = \ord d + \dual{\mu
    -w(\mu)}{\rho\spcheck}$ for all $0 \ne d \in \cD^\mu$. In
  particular, if $0 \ne f \in \cO^\eta$, then $\ord \Fo(f) =
  \dual{\eta}{2\rho\spcheck}$.
\item[(4)] Let $\{f_i\}_{1\le i \le n} \subset \cO^\eta$ and
  $\{g_i\}_{1\le i \le n} \subset \cO^{\eta^*}\cong
  (\cO^\eta)^*$ be dual bases such that (up to a constant)
  $\sum_i f_i \otimes g_i$ is the unique $G$-invariant element
  in $\cO^\eta\otimes \cO^{\eta^*}$.  Then the $(G \times
  H)$-invariant operator $P_{\eta} = \sum_{i=1}^n f_i \Fo(g_i)
  \in U(\wh)$ is given by
\begin{equation}\label{eq-prop22}
P_\eta = c_\eta \prod_{\alphacheck \in
    \Delta_{+}\spcheck} \prod_{i=1}^{\dual{\eta}{\alphacheck}}
  (\alphacheck + \dual{\alphacheck}{\rho}-i),
\end{equation}
for some $c_\eta\in \C\smallsetminus\{0\}$.  Moreover, $\ord
P_\eta =\dual{\eta}{2\rho\spcheck}$.  \vskip 3pt
\item[(5)] \label{lastitem} $U(\wh) = \sum_{w\in
    W}U(\wh)F_w(P_\eta)$. \qed
\end{enumerate}
\end{prop}


\section{The structure of $\cD(G/U)$.}
\label{structure-section}

In \cite{BBP}, the authors use Theorem~\ref{thm21} to prove
that $\cD$ is a noetherian ring of finite injective dimension.
In this section we investigate other consequences of that
result to the structure of $\cD$.  The notation from the last
section will be retained; in particular, $G$ is a connected,
simply connected reductive algebraic group over $\mathbb C$,
with basic affine space $\mbX=G/U$ and $\cD=\cD(\mbX)$.

We begin with an easy application of the faithful flatness of
the ring $R=\bigoplus_{w,j}\cD(\mbX_j)_{w}$ defined in
Corollary~\ref{cor24}.

\begin{prop}
\label{cor26}
The ring $\cD$ is simple and $\cO$ is a simple left
$\cD$-module.
\end{prop}

\begin{proof}
  Let $J$ be a non zero ideal of $\cD$. As in the proof of
  Corollary~\ref{cor24}, $\cD(\mbX_j)$ is a noetherian
  localization of $\cD$ and so, by
  \cite[Proposition~2.1.16(vi)]{MR}, each
  $\cD(\mbX_j)_{w}\otimes_\cD J \cong \cD(\mbX_j)\otimes_\cD
  J^{w}$ is an ideal of $\cD(\mbX_j)$ and it is nonzero since
  $\cD$ is a domain. But $\mbX_j$ is a smooth affine variety
  and so $\cD(\mbX_j)$ is a simple ring. Thus
  $\cD(\mbX_j)_{w}\otimes_\cD J = \cD(\mbX_j)_{w}$ for all
  $w\in W$ and $1\leq j\leq k$.  This means that the module $M
  = \cD/J$ satisfies $R \otimes_\cD M = 0$, hence $M=0$ by
  Corollary~\ref{cor24}(2).  In other words, $J = \cD$.
  
  If $\cO$ is not a simple $\cD$-module, pick a proper factor
  module $\cO/K$ and note that $K$ is then an ideal of $\cO$.
  But now the annihilator $\mathrm{ann}_{\cD}(\cO/K)$ of
  $\cO/K$ as a $\cD$-module is an ideal of $\cD$ that contains
  $K$.  This contradicts the simplicity of $\cD$.
\end{proof}

We now turn to the homological properties of $\cD$.  Two
conditions that have proved very useful in applying homological
techniques (see for example \cite{Bj1} or \cite{LS0}) are the
Auslander and Cohen-Macaulay conditions.  These are defined as
follows.  A noetherian algebra $A$ of finite injective
dimension is called \emph{Auslander-Gorenstein} if, for any
finitely generated (left) $A$-module $M$ and any right
submodule $N\subseteq \Ext_A^i(M,A)$, one has
$\Ext^j_A(N,A)=0$, for $j<i$.  The \emph{grade} of $M$ is
$\grade_A(M)= \inf \{j : \Ext_A^j(M,A) \ne 0\}$ (with the
convention that $\grade_A(0) = +\infty$).  The
\emph{Gelfand-Kirillov dimension} of $M$ will be denoted
$\GK_AM$.  We say that the algebra $A$ is \emph{Cohen-Macaulay}
if
$$
\GK_A M + \grade_A(M) = \GK A \ \; \text{for all $M \ne 0$}.
$$

If $\mbZ$ is a smooth affine variety, then $\cD(\mbZ)$ is
Auslander-Gorenstein and Cohen-Macaulay with $\GK \cD(\mbZ) = 2
\dim \mbZ$ (see \cite[Chapter~2, Section~7]{Bj1}).  This
applies, of course, when $\mbZ=\mbX_j$ for $1\leq j\leq k$ and
so $R = \bigoplus_{w,j}\cD(\mbX_j)_{w}$ also satisfies these
properties.  As we show in Theorem~\ref{thm29}, these
properties descend to $\cD=\DX$.

If $M$ is a (finitely generated) $\cD$-module, write
$$M_j^w = \cD(\mbX_j) \otimes_\cD M^w, \qquad\text{for}\quad 1
\le j \le k\ \text{and}\ w \in W.$$

\begin{lem}
\label{lem28}   
If $M$ is a finitely generated left $\cD$-module then
$$
\GK_\cD M = \max\{\GK_{\cD(\mbX_j)} M_j^w : 1 \le j \le k, w
\in W\}.
$$
\end{lem}

\begin{proof}
  By definition, $R \otimes_\cD M = \bigoplus_{j=1}^k
  \bigoplus_{w \in W} M_j^w$ and so
  $$
  \GK_R R \otimes_\cD M = \max\{\GK_{\cD(\mbX_j)} M_j^w : 1
  \le j \le k, w \in W\}.
  $$
  By faithfully flatness (Corollary~\ref{cor24}), the
  natural map $M \to R \otimes_\cD M$ is injective and it
  follows that $\GK_\cD M \le \GK_R R \otimes_\cD M$.
  Conversely, since $M_j^w$ is the localization of $M^w$ at the
  Ore subset $\{f_j^s : s \in \N\}$, it follows from
  \cite[Theorem~3.2]{Lo} that $\GK_{\cD(\mbX_j)} M_j^w \le
  \GK_\cD M^w$. Since $\GK_\cD M^w = \GK_\cD M$, this implies
  that $\GK_R R \otimes_\cD M \le \GK_\cD M$ and the lemma is
  proved.
\end{proof}

\begin{thm}
\label{thm29}
The algebra $\cD$ is Auslander-Gorenstein and Cohen-Macaulay.
\end{thm}

\begin{proof}
  Let $M$ be a finitely generated (left) $\cD$-module and $N$
  is a (right) submodule of $\Ext_\cD^p(M,\cD)$.  As $R$ is a
  flat $\cD$-module, $$R \otimes_\cD \Ext_\cD^i(N,\cD) \cong
  \Ext_R^i(N \otimes_\cD R,R)$$
  and $N\otimes_\cD R$ is a
  submodule of $\Ext_R^p( R \otimes_\cD M,R)$.  Since $R$ is
  Auslander-Gorenstein, this implies that $R \otimes_\cD
  \Ext_\cD^q(N,\cD)=0$, for any $q<p$. Since $R$ is a faithful
  $\cD$-module, it follows that $\Ext_\cD^q(N,\cD)=0$. Thus
  $\cD$ is Auslander-Gorenstein.
  
  From $\Ext_R^i(R \otimes_\cD M, R) \cong
  \Ext_\cD^i(M,\cD)\otimes_\cD R $ and the faithful flatness of
  $R_{\cD}$ we get that $ \Ext_R^i(R \otimes_\cD M, R) = 0 \iff
  \Ext_\cD^i(M,\cD) =0$. Thus, $\grade_{R}(R \otimes_\cD M) =
  \grade_\cD(M)$. But,
  $$
  \Ext_R^p(R \otimes_\cD M, R)\ \cong\ \bigoplus_{j,w}\,
  \Ext_{\cD(\mbX_j)}^p(M_j^w, \cD(\mbX_j))
  $$
  and so $\grade_\cD(M)=\grade_{R}(R \otimes_\cD M) = \min\{
  \grade_{\cD(\mbX_j)}(M_j^w) : 1 \le j \le k, w \in W\}$.
  
  Each $\mbX_j$ is a smooth affine variety of dimension $\dim
  \mbX$ and so each $\cD(\mbX_j)$ is Cohen-Macaulay.  It
  therefore follows from Lemma~\ref{lem28} that
\begin{equation*}
\begin{split}
  \GK_\cD M &= \max\{\GK_{\cD(\mbX_j)} M_j^w : 1 \le j \le
  k, w \in W, M_j^w \ne 0\} \\
  &=\max\{2\dim \mbX - \grade_{\cD(\mbX_j)}(M_j^w) : 1 \le j
  \le
  k, w \in W, M_j^w \ne 0\}  \\
  &= 2\dim \mbX - \min\{
  \grade_{\cD(\mbX_j)}(M_j^w) : 1 \le j \le k, w \in W\} \\
  &= 2\dim \mbX - \grade_\cD(M),
\end{split}
\end{equation*}
as required.
\end{proof}

\begin{rem}\label{YZ-remark} The last result should be compared
  with \cite[Corollary~0.3]{YZ} which shows that a simple
  noetherian $\C$-algebra with an affine commutative associated
  graded ring is automatically Auslander-Gorenstein and
  Cohen-Macaulay. It is not clear whether that result applies
  to $\cD$, since we do not have a good description of the
  associated graded ring of $\cD$.
\end{rem}

When $\mbZ$ is a smooth affine variety, $\cD(\mbZ)$ is a
finitely generated $\mathbb C$-algebra, simply because it is
generated by $\cO(\mbZ)$ and $\mathrm{Der}_{\mathbb C}(\mbZ)$.
When $\mbZ$ is singular, it can easily happen that $\cD(\mbZ)$
is not affine (see, for example \cite{BGG}). However, as we
will show in Theorem~\ref{thm36}, the $\C$-algebra $\cD$ is
finitely generated. We begin with some lemmas.

\begin{lem}
\label{lem23}
Let $\gamma \in \Lambda^+$. There exists a nondegenerate
pairing of $G$-modules
$$
\tau: \Fo(\cO^\gamma) \times \cO^{\gamma^*} \lto \C, \quad
(\Fo(f),g) \mapsto \Fo(f)(g).
$$
\end{lem}

\begin{proof}
  By Proposition~\ref{prop22} and \eqref{D-mu},
  $\Fo(\cO^\gamma)(\cO^{\gamma^*}) \subseteq \cO^{w_0(\gamma) +
    \gamma^*} = \cO^0 = \C$. Also,
  $$
  (a.\Fo(f))(a.g) = a.\Fo(f)(a^{-1}.(a.g)) = a.\Fo(f)(g)
  =\Fo(f)(g),
  $$
  for all $a \in G$. Therefore, $\tau$ is a well defined
  pairing of $G$-modules and it induces a $G$-linear map
  $\Fo(\cO^\gamma) \to (\cO^{\gamma^*})^*$.  Since
  $\Fo(\cO^\gamma) \cong (\cO^{\gamma^*})^* \cong V(\gamma)$,
  it suffices to show that $\tau$ is non zero in order to show
  that it is non-degenerate.  In the notation of
  Proposition~\ref{prop22}(4), we will show that $\Fo(f_i)(g_i)
  \ne 0$ for some $i$.
  
  By Proposition~\ref{prop22}(2),
  $$
  \Fo(\alphac) \ = \ w_0(\alphac) + \dual{w_0(\alphac) -
    \alphac}{\rho} \ = \ w_0(\alphac) - 2\dual{\alphac}{\rho}$$
  and so Proposition~\ref{prop22}(4) implies that
  $$
  {\Fo(P_\gamma) \ = \ c_\gamma \prod_{\alphacheck \in
      \Delta_{+}\spcheck}
    \prod_{i=1}^{\dual{\gamma}{\alphacheck}} (w_0(\alphacheck)
    - \dual{\alphacheck}{\rho}-i)}.
  $$
  Since $w_0(\alphac) \in \wh$ is a vector field on $\mbX$,
  we have $w_0(\alphac)(1) = 0$ and so
  $$
  {\Fo(P_\gamma)(1) \ = \ c_\gamma \prod_{\alphacheck \in
      \Delta_{+}\spcheck} (-1)^{\dual{\gamma}{\alphacheck}}
    \prod_{i=1}^{\dual{\gamma}{\alphacheck}}(\dual{\alphac}{\rho}
    +i) \ \ne\ 0, }
  $$
  since $\dual{\alphac}{\rho} > 0$ for all $\alphac$. By
  Proposition~\ref{prop22}(4), $\Fo(P_\gamma)(1)=\sum_i
  \Fo(f_i)(g_i) $ and so $\Fo(f_i)(g_i) \ne 0$ for some $1 \le
  i \le n$.
\end{proof}

\begin{lem}
\label{lem34}
The $\cO(\mbX)$-module map $ \mathsf{m} : \cO(\mbX) \otimes
\Fo(\cO(\mbX)) \to \cD(\mbX)$ given by $ g\otimes\Fo(f) \mapsto
g\Fo(f)$ is injective.
\end{lem}

\begin{proof}
  Let $t \in \Ker(\mathsf{m})$ and write $t=\sum_{i,j} h_{ij}
  \otimes \Fo(f_{ij})$ where $\{f_{ij}\}_j$ is a basis of
  $\cO^{\mu_i}$, for some $\mu_1,\dots, \mu_s \in \Lambda^+$,
  and $h_{ij} \in \cO$ for all $i,j$. We may assume that $s$ is
  minimal; thus, for each $1\leq i\leq s$, some $h_{ij}\not=0$.
  Partially order $\Lambda$ by $\omega^*\ges \lambda^*$ if
  $\omega^*-\lambda^* \in \Lambda^+$ and assume that $\mu_1^*$
  is minimal among the $\mu_i^*$'s for this ordering. By
  Lemma~\ref{lem23}, for each $i$ there exists a basis
  $\{g_{ij}\}_j$ of $\cO^{\mu_i^*}$ such that
  $\Fo(f_{ij})(g_{ik}) = \delta_{jk}$. When $i > 1$ we have
  $\Fo(f_{ij})(g_{1k}) \subseteq \cO^{\mu_1^* - \mu_i^*} = 0$,
  since $\mu_1^* - \mu_i^* \notin \Lambda^+$. Therefore, for
  each $k$, we have
  $$
  {0 =\mathsf{m}(t)(g_{1k})= \sum_{i,j}
    h_{ij}\Fo(f_{ij})(g_{1k}) = \sum_j
    h_{1j}\Fo(f_{1j})(g_{1k}) = h_{1k} },
  $$
  contradicting the minimality of $s$.
\end{proof}

Define a finitely generated subalgebra of $\DX$ by
$$
\calS= \C \langle \cO, F_{w_0}(\cO) \rangle = \C \langle
\cO^{\varpi_i}, \Fo(\cO^{\varpi_j}) \, ; \, 1 \le i,j \le \ell
\rangle,
$$
where $\ell=\text{rank}(\g)$. The elements of $\cO$ act
locally nilpotently on $\cD$, and therefore on $\calS$. Thus
$C=\cO \smallsetminus \{0\}$ is an Ore subset of $\calS$. Let
$\K= C^{-1}\cO$ denote the field of fractions of $\cO$. Recall
that $C^{-1} \DX = \cD(\K)$ is the ring of differential
operators on $\K$ and that $\cD_r(\mbX) = \cD_r(\K) \cap
\cD(\mbX)$ for all $r \in \N$ (see, for example,
\cite[Theorem~15.5.5]{MR}).

\begin{lem}
\label{lem35}
We have $C^{-1}\calS = \cD(\K)$. In particular, for any finite
dimensional subspace $E \subset \DX$, there exists $0 \ne f \in
\OX$ such that $fE \subset \calS$.
\end{lem}

\begin{proof} The aim of the proof is to apply \cite[Lemma~8]{LS1}.
  
  Applying the exact functor $\K \otimes_\cO {-}$ to the
  injective map $\mathsf{m}$ of Lemma~\ref{lem34} yields the
  $\K$-linear injection $\mathsf{m} : \K \otimes\Fo(\cO) \ito
  A=C^{-1}\calS$. Since $\Fo(\cO)$ is a commutative algebra of
  dimension $n = \dim \mbX = \trdeg_\C \K$, we may pick
  $u_1,\dots,u_n \in \Fo(\cO)$ algebraically independent over
  $\C$ and set $P = \C[u_1,\dots,u_n]$. For any $q \in \N$,
  denote by $P_q$ the subspace of polynomials of degree at most
  $q$ and define
  $$
  \Theta_q = \{ d \in \K \otimes P : \ord \mathsf{m}(d) \le
  q\} \quad\text{and}\quad k = \max\{\ord \mathsf{m}(u_i) : 1
  \le i \le n\}.
  $$
  Observe that $\{\Theta_q\}_q$ and $\{\K \otimes P_q\}_q$
  are two increasing filtrations on the $\K$-vector space $\K
  \otimes P$ and that the map $q \mapsto \dim_\K (\K \otimes
  P_q)$ is a polynomial function of degree $n$. Furthermore,
  since $\mathsf{m}(u_iu_j)=\mathsf{m}(u_i)\mathsf{m}(u_j)$, we
  have $\K \otimes P_q \subseteq \Theta_{kq}$ for all $q$.
  Hence, $\dim_\K \Theta_q \ge p(q)$ for some polynomial $p$ of
  degree $n$. If $A_q = \cD_q(\K) \cap A$, then
  $\mathsf{m}(\Theta_q) \subseteq A_q$ and so $\dim_\K A_q \ge
  p(q)$. It follows that
  $$
  \limsup_{q \to \infty} \{\log_q(\dim_\K A_q/A_{q-1})\} \ge
  n-1.
  $$
  
  The hypotheses of \cite[Lemma~8]{LS1} are now satisfied by
  the pair $A \subseteq \cD(\K)$ and, by that result, $\cD(\K)=
  A = C^{-1}\calS$.  The final assertion of the lemma follows
  by clearing denominators.
\end{proof}

We can now describe a generating set for the $\C$-algebra
$\cD$, for which we recall the definition of $\widehat{\h}$
from Notation~\ref{h-hat}.

\begin{thm}
\label{thm36}
As a $\C$-algebra, $\cD$ is generated by $\widehat{\h}$ and the
$F_w(\cO)$, for $w\in W$.
\end{thm}

\begin{proof}    Set $\calB = \C \langle\, \widehat{\h}, F_w(\cO)
  \, ; \, w \in W \rangle$; thus $\calB$ is a $G$-submodule of
  $\cD$ containing both $\calS$ and $U(\wh)$.  Moreover, as
  $F_w(U(\wh)) = U(\wh)$ (see Proposition~\ref{prop22}(2)),
  $F_w(\calB) = \calB$ for all $w \in W$.  As $\cD$ is a
  locally finite $G$-module, it suffices to show that $E
  \subset \calB$ for all finite dimensional $G$-submodules $E$
  of $\cD$. For such a module $E$, set $ L = \{b \in \calB : b
  E \subset \calB \}$ and $I_w(E) = \{g \in \cO : F_w(g) E
  \subset \calB\}$ for $w\in W.  $ We aim to show that the left
  ideal $L$ of $\calB$ contains $1$.
  
  Clearly, $I_w(E)$ is an ideal of $\cO$. It is also a
  $G$-submodule since $$F_w(a.g)E = F_w(a.g)a.E = a.(F_w(g)E)
  \subset a.\calB = \calB$$
  for all $a \in G$ and $ g \in
  I_w(E)$.  Since $\calS \subseteq \calB$, Lemma~\ref{lem35}
  implies that $I_1(E)\not=0$.  For each $w$, $F_{w^{-1}}(E)$
  is a $G$-submodule of $\cD$ (isomorphic to $E$) and $$g \in
  I_w(E) \iff g F_{w^{-1}}(E) \subset \calB \iff g \in
  I_{1}(F_{w^{-1}}(E)).$$
  Thus, $I_w(E) \ne 0$ for all $w \in
  W$. Since $\cO$ is a domain, it follows that $I = \bigcap_{w
    \in W} I_w(E)$ is a non-zero $G$-submodule of $\cO$. Now,
  $\cO = \bigoplus_{\lambda \in \Lambda^+} \cO^\lambda$ and so
  $\cO^\gamma \subset I$ for some $\gamma \in \Lambda^+$. By
  Proposition~\ref{prop22}(4), $F_w(P_{\gammas}) = \sum_i
  F_w(f_i) F_{w w_0}(g_i)$ for some $f_i \in \cO^{\gammas}$ and
  $ g_i \in \cO^\gamma$.  By the definition of $I$ and the
  choice of $\gamma$, we have $F_{w w_0}(g_i)E \subset \calB$;
  that is, $F_{w w_0}(g_i) \in L$, for all $w \in W$. Since
  $F_w(f_i) \in F_w(\cO) \subset \calB$ we obtain that
  $F_w(P_{\gammas}) \in L$ for all $w \in W$. Finally, as
  $U(\wh)\subset \calB$, Proposition~\ref{prop22}(5) implies
  that $1 \in L$ and hence that $E\subset \calB$.
\end{proof}

We conjecture that $\widehat{\h}$ is unnecessary in the last
theorem; i.e., we conjecture that
$$
\cD=\C \langle F_w(\cO) : w\in W \rangle.
$$
Using Example~\ref{sl-3} the authors can prove this for
$\g=\mathfrak{sl}(3)$; indeed we can even show that $\cD=\C
\langle \cO, F_{w_0}(\cO) \rangle$ in this case.  However, the
argument heavily uses facts about $\mathfrak{so}(8)$ and so it
may not be a good guide to the general case.


\section{The $\DX$-module $\coh^*(\mbX,\cO_\mbX)$}
\label{sec5}

As before, set $\mbX=G/U$ with associated rings $\cO=\cO(\mbX)$
and $\cD=\cD(\mbX)$.  In this section we give a complete
description of $ \coh^*(\mbX,\cO_\mbX) = \bigoplus_{i \in \N}
\coh^i(\mbX,\sOX)$ as a $\cD$-module.  Specifically, we will
show that $\coh^*(\mbX,\cO_\mbX)$ is simply the direct sum of
the twists $\cO^w=\cO^{F_w}$ of $\cO$ by $w\in W$.

The first cohomology group to consider is $\coh^0(\mbX,\sOX) =
\cO$. As a $\cD$-module, $\cO \cong \cD/I$ where $I = \{d \in
\cD : d(1) = 0\}$ is a maximal left ideal of $\cD$
(Proposition~\ref{cor26}).  Clearly, $I \supseteq \cD \g + \cD
\widehat{\h}$, in the terminology of Notation~\ref{h-hat}.  The
first main result of this section, Theorem~\ref{thm54}, shows
that this is actually an equality and, further, that
$\cD/\cD\g\cong \bigoplus_{w\in W} \cO^w$.

Before proving this theorem, we need some preliminary notation
and lemmas.  As in \eqref{cover-equ}, we cover $\mbX$ by affine
open subsets $\mbX_j = \{x \in \overline{\mbX} : f_j(x) \ne
0\}$ and let $\calC_j = \{f_j^s : s \in \N\}$, for $1 \le j \le
k$, denote the associated Ore subsets in $\cD$.  If $M$ is a
left $\cD$-module, the kernel of the localization map $M \to
M_{f_j} = \cD(\mbX_j) \otimes_\cD M$ is
$$
T_j(M) = \{v \in M : f_j^s v = 0 \ \text{for some } s>0\}.
$$

\begin{lem}
\label{g-fixed} 
For all $w\in W$ and $x\in \g$, one has $F_w(x)= x$.
\end{lem}

\begin{proof} 
  By Proposition~\ref{prop22}(1), $F_w$ is a $G$-linear
  automorphism of $\cD$ and hence is $\g$-linear, where $\g$
  acts by the adjoint action.  Therefore, for any $\theta\in
  \cD$ and $x \in \g$,
  $$
  [F_w(x),F_w(\theta)]=F_w([x,\theta]) =[x,F_w(\theta)].
  $$
  Thus $[F_w(x),y]=[x,y]$ for all $y \in \cD$.  If $y\in
  \g$, then $\g$-linearity also implies that $[F_w(x),
  y]=F_w([x,y])$ and hence that $[x,y]=[F_w(x),y]=F_w([x,y])$.
  As $\g$ is semisimple, $[\g,\g]=\g$ and so this implies
  $F_w(z) = z$ for all $z \in \g$.
\end{proof}

The next lemma shows why we should expect all the $\cO^w$ to
appear in a decomposition of $\cD/\cD\g$.

\begin{lem}
\label{lem51}
{\rm (1)} Each $\cD$-module $\cO^w$ is a factor of $\calM = \cD
/ \cD \g$.

\noindent{\rm (2)} Set $\calN = I/\cD \g$. Then, $\calN =
T_j(\calM)$ for $1 \le j \le k$ and $\cD/ I$ is the unique
factor of $\calM$ isomorphic to $\cO$ as a $\cD$-module.
\end{lem}

\begin{proof}  (1) Since the action of $d \in \cD$ on $\cO^w$
  is given by $d\cdot f = F_w(d)(f)$ for $f \in \cO$, we have
  $\cO^w \cong \cD / F_w^{-1}(I)$. As we noted above,
  $I\supseteq \cD\g +\cD \widehat{\h}$ and so, by
  Lemma~\ref{g-fixed}, $F_w^{-1}(I) \supseteq \cD \g$.

  (2) Let $T_x\mbX$ denote the tangent space of $\mbX$ at $x
  \in \mbX$. Observe that, if $e=U/U \in \mbX$, then
  $\mathrm{Stab}_G(e)=U$ and so $\g / \fu$ identifies naturally
  with $T_{e}\mbX$. Since $\mbX$ is a homogeneous space, the
  map $\imath : \g \to \Der_\C\cO$ induces an isomorphism $\g /
  \fu \cong T_e\mbX \cong T_x\mbX$.  As $\mbX_j \subseteq \mbX$
  is affine, it follows that $\Der \cO(\mbX_j) = \cO(\mbX_j)
  \g$ for each $j=1,\dots,k$.  Furthermore, since $\mbX_j$ is
  smooth, the $\cD(\mbX_j)$-module $\cO(\mbX_j)$ is simple and
  isomorphic to $\cD(\mbX_j)/\cD(\mbX_j) \Der \cO(\mbX_j)$.
  Since $\cO = \calM/\calN$, this implies that
  \begin{equation}
  \label{lem511}
  \cO(\mbX_j) = \calC_j^{-1} (\calM/\calN) \cong \cD(\mbX_j)
  / \cD(\mbX_j) \g = \calC_j^{-1} \calM
  \end{equation}
  for all $j$. Therefore $\calC_j^{-1} \calN = 0$ and so $\calN
  \subseteq T_j(\calM)$. The equality $\calN = T_j(\calM)$ then
  follows from $T_j(\cO)=0$.
  
  Now suppose that $\cO \cong \cD/ L$ for some maximal left
  ideal $L \supseteq \cD \g$. Then $\calC_j^{-1} L \supseteq
  \cD(\mbX_j) \g$ and $\calC_j^{-1}\cO = \cO(\mbX_j) \cong
  \cD(\mbX_j)/ \calC_j^{-1} L$ for all $j$. By \eqref{lem511},
  the left ideal $\cD(\mbX_j) \g$ is maximal and so
  $\calC_j^{-1} L = \cD(\mbX_j) \g$ for all $j$.  Therefore
  $L/\cD \g \subseteq T_j(\calM) = \calN = I/\cD \g$, which
  implies that $L \subseteq I$. Hence $L=I$.
\end{proof}

As in \cite[Th\'eor\`eme~2, Section~VI.1.5]{Bou},
$\Sigma\spcheck =\{\alpha_1\spcheck, \dots,
\alpha_\ell\spcheck\}$ defines a dominant chamber
$$
C(\Sigma\spcheck) = \bigl\{y \in \wh_\R : \dual{y}{\alpha_i}
> 0 \, \; \text{for all $i =1,\dots,\ell$}\bigr\}
$$
in the root system $\Deltacheck = \jmath(\Deltacheck)
\subset \wh_\R = \bigoplus_{i=1}^\ell \R \alpha_i\spcheck
\subset \wh$.

\begin{lem}
\label{lem52}
If $y \in C(\Sigma\spcheck)$ and $w\in W\smallsetminus\{ 1\}$,
then $F_w(y) \notin I$.
\end{lem}

\begin{proof} Using  \cite[Ch.~VI, Section~1.6,
  Corollaire de la Proposition~18]{Bou}, we have $0 \ne y -
  w(y) = \sum_{i=1}^\ell n_i \alpha_i\spcheck$ with $n_i \in
  \R_+$.  Hence, Proposition~\ref{prop22}(2) implies that
  $$F_w(y) - w(y)\ =\ \dual{w(y)-y}{\rho}\ =\ -\textstyle\sum_i
  n_i\ \in\ \C\smallsetminus\{0\}.$$
  Thus $F_w(y) - w(y) \notin
  I$.  Since $w(y)\in\wh \subset I$, this implies that
  $F_w(y)\notin I$.
\end{proof}

\begin{lem}
\label{lem53}
As $\cD$-modules, $\cO^w \cong \cO^{w'}$ if and only if $
w=w'$.
\end{lem}

\begin{proof}
  Since $(\cO^w)^{v} \cong \cO^{wv}$ it suffices to prove the
  result when $w'=1$. So, assume that $\cO^w \cong \cO$ for
  some $w \ne 1$ and set $\calN_w = F_{w}^{-1}(I)/ \cD \g$;
  thus $\calN_1=\calN$.  Then $\cO^w = \cD / F_w^{-1}(I) =
  \calM / \calN_w$ and Lemma~\ref{lem51} implies that $\calN_w
  = \calN_1$; equivalently, $ F_{w^{-1}}(I) = I$.  Now pick $y
  \in C(\Sigma\spcheck)$.  Since $y\in I$, this implies that
  $F_{w^{-1}}(y)\in I$, contradicting Lemma~\ref{lem52}.
\end{proof}

The following theorem gives a precise description of the
$\cD$-modules $\cO$ and $\calM$; it shows in particular that
$\calM$ is a multiplicity free, semisimple module of length
$|W|$.

\begin{thm}
\label{thm54}  Write  $\OX \cong \DX / I$ for $I=\{d \in
\DX : d(1) = 0\}$ and define $\calM = \DX/ \DX \g$. Then:
\begin{enumerate}
\item[(1)] $\calM \cong \bigoplus_{w \in W} \OXw$ as a
  $\DX$-module;
\item[(2)] $I = \DX \g + \DX \wh = \DX\g + \DX y$ for all $y
  \in C(\Sigma\spcheck)$.
\end{enumerate}
\end{thm}

\begin{proof}
  (1) Set $\calN_w = F_w^{-1}(I)/\cD \g$, for $w \in W$. By
  Lemmas~\ref{lem51}(2) and~\ref{lem53}, the $\cD$-modules
  $\cO^w\cong \calM/\calN_w$ are nonisomorphic, so the natural
  map $\calM \to \bigoplus_{w \in W} \cO^w$ is surjective with
  kernel $N=\bigcap_{w \in W} \calN_w$. It therefore remains to
  prove that $N=0$.  By Lemma~\ref{lem51}, if $x \in I$, then
  there exists $t \in \N$ such that $f_j^t x \in \cD \g$ for
  $j=1,\dots,k$. Therefore $F_{w^{-1}}(f_j)^t F_{w^{-1}}(x) \in
  \cD \g$ for all $j$; equivalently, each element
  $[F_{w^{-1}}(x)+\cD\g ] \in \calN_w$ is torsion for
  $F_{w^{-1}}(\calC_j)$.  Consequently, if $v \in N$, then
  there exists $s \in \N$ such that $F_w(f_j)^sv = 0$ for all
  $1 \le j \le k$ and $w \in W$. In other words,
  $N^w_{f_j}=(0)$, for all such $j$ and $w$. By faithful
  flatness (Corollary~\ref{cor24}) this implies that $N=0$.
  
  (2) It suffices to prove that $I = \cD\g + \cD y$ for $y \in
  C(\Sigma\spcheck)$. Set $\eI= (\cD\g + \cD y)/\cD \g$.  By
  part~(1) and Lemma~\ref{lem51}, the $\calN_w$, $w \in W$, are
  the only maximal submodules of $\calM$. Therefore, if $\eI
  \subsetneq \calN_1$ we must have $\eI \subseteq \calN_w$ for
  some $w \ne 1$.  This implies that $F_w(\cD \g + \cD y)
  \subseteq I$, hence $F_w(y) - w(y) \in I$, in contradiction
  with Lemma~\ref{lem52}.  Therefore $\calN_1=\eI$ and $I =
  \cD\g + \cD y$.
\end{proof}

We will prove in Theorem~\ref{thm57} that $\coh^*(\mbX,\sOX)$
is isomorphic to $\calM$ as a $\cD$-module. In order to prove
this, we need to recall some results on the cohomology of line
bundles over the flag variety $\mbB=G/B$.

We begin with some general remarks. Let $\mbZ$ be a smooth
$G$-variety and write $\tau : \g \to \Der \cO_\mbZ$ for the
differential of the $G$-action. Let $\eF$ be a $G$-equivariant
$\cO_\mbZ$-module as defined, for example, in
\cite[Section~{4.4}]{Ka}.  This implies, in particular, that
$\eF$ is a compatible $(\g,\cO_\mbZ)$-module in the sense that
for any open subset $\Omega \subseteq \mbZ$, one has $\xi.(fv)
= \tau(\xi)(f)v + f(\xi.v)$ for all $\xi \in \g, f \in
\cO_\mbZ(\Omega)$ and $ v \in \eF(\Omega)$.  In this setting,
the cohomology group $\coh^{i}(\mbZ,\eF)$ inherits a structure
of compatible $(G,\cO(\mbZ))$-module,
cf.~\cite[Theorem~11.6]{Ke2}.  If $\eF$ is a coherent
$G$-equivariant $\cD_\mbZ$-module, then it follows from
\cite[Sections~{4.10} and {4.11}]{Ka} that $\coh^{i}(\mbZ,\eF)$
is endowed with a $\cD(\mbZ)$-module structure such that
$\tau(\xi)v = \frac{d}{dt}_{\mid t=0} (e^{t \xi}.v)$ for all $v
\in \coh^{i}(\mbZ,\eF)$ and $ \xi \in \g$. We will apply these
observations in two cases: one is when $\mbZ = \mbX = G/U$ and
$\eF=\cO_\mbX$ is a $G$-equivariant $\cD_\mbX$-module under
left translation; the other is described next.
  
Let $\pi : G \sto \mbB=G/B$, $\phi : G \sto \mbX$ and $\vphi :
\mbX \sto \mbB$ be the natural projections, thus $\pi = \vphi
\circ \phi$.  For each $\lambda \in \Lambda$, the one
dimensional $H$-module $\C_{-\lambda}$ can be viewed as a
$B$-module with trivial action of $U$.  As in
\cite[pp.333-335]{Ke2}, define the $G$-equivariant
$\cO_\mbB$-module $\L{\lambda}$ to be the sections of the line
bundle $G \times^B \C_{-\lambda}$. Since $G\times^B
\C_{-\lambda} \cong \mbX \times^H \C_{-\lambda}$, one has
\begin{equation}\label{L-sheaf}
\Gamma(\Omega,\L{\lambda}) = \{f : \vphi^{-1}\Omega \to \C
: f(\bar{g}h) = \lambda(h)f(\bar{g}) \; \, \text{for $h \in
  H$, $\bar{g} \in \vphi^{-1}\Omega$}\}.
\end{equation}
where $\Omega \subseteq \mbB$ is any open subset.

Recall from \eqref{O-decomp} that the decomposition $\cO =\OX=
\bigoplus_{\gamma \in \Lambda^+} \cO^\gamma$ is induced by the
twisted action of $H$ on $\mbX$.  Hence, $h \in H$ acts on $f
\in \Gamma(\Omega,\L{\lambda})$ via
\begin{equation}\label{rem541}
(r_h.f)(\bar{g})= f(\bar{g}w_0(h^{-1})) =
\lambda(w_0(h)^{-1})f(\bar{g}) = \lambda^*(h)f(\bar{g}).
\end{equation}
Passing to \v{C}ech cohomology, $H$ therefore acts on
$\coh^i(\mbB, \L{\lambda})$ with weight $\lambda^*$.

The cohomology groups of the line bundle $\L{\lambda}$ are
described as $G$-modules by the Borel-Weil-Bott theorem, that
we now recall (see \cite[Corollaries~II.5.5 and II.5.6]{Ja}, up
to a switch from $B$ to the opposite Borel). The length of $w
\in W$ is denoted by $\ell(w)$ and we will write
$$
W(i)=\{w\in W : \ell(w) = i\}.
$$
The ``dot action'' of $w \in W$ on $\xi \in \fh^*$ is
defined by $w\cdot \xi = w(\xi + \rho) - \rho$.

\begin{thm}
\label{thm55} \emph{(Borel-Weil-Bott)}
The $G$-module $\coh^i(\mbB,\L{\lambda})$ is isomorphic to
$$
\qquad
\begin{cases}
  V(\mu)^* &\text{if $\, \exists \, (w,\mu) \in W(i) \times
    \Lambda^+$
    such that $\lambda = w\cdot \mu$;} \\
  0 &\text{otherwise.}  \hfill {}
\end{cases} \quad\begin{array}{c}{\strut}\\
    \textstyle\qed\end{array} 
  $$
\end{thm}

The following standard proposition reduces the computation of
the $G$-module $\coh^*(\mbX,\sOX)$ to the Borel-Weil-Bott
theorem. We include a proof since we could not find an
appropriate reference.

\begin{prop}
\label{prop56}
The morphism $\vphi : \mbX \to \mbB$ is affine and $\vphi_*\sOX
\cong \bigoplus_{\nu \in \Lambda} \L{\nu}$ as
$(\g,\cO_\mbB)$-modules. In particular, for each $i \in \N$,
there is a $G$-module isomorphism
$$
\coh^i(\mbX,\sOX) \cong {\bigoplus_{\nu \in \Lambda}}
\coh^i(\mbB,\L{\nu}).
$$
\end{prop}

\begin{proof}
  By the Bruhat decomposition, $\mbB$ is covered by the affine
  open subsets $\Omega_w= \pi(\dot{w}U^-B)$ where $\dot{w} \in
  N_G(H)$ is a representative of $w \in W$ and $U^-$ is the
  opposite maximal unipotent subgroup of $G$, see
  \cite[(II.1.10)]{Ja}.  As $\dot{w}U^-B \cong U^- \times H
  \times U$, the subset $\vphi^{-1} \Omega_w =
  \phi(\dot{w}U^-B)$ is affine and isomorphic to $U^- \times H$
  as an $H$-variety. Thus $\vphi$ is an affine morphism and
  $\coh^i(\mbX,\sOX) \cong \coh^i(\mbB,\vphi_*\sOX)$ by
  \cite[III, Exercise~4.1]{Ha}.
  
  The affine algebra $\Gamma(\Omega_w, \vphi_*\sOX) =
  \sOX(\vphi^{-1}\Omega_w)$ is endowed with a regular action of
  $H$.  Hence, $\Gamma(\Omega_w, \vphi_*\sOX)$ decomposes as
  $\bigoplus_{\nu \in \Lambda} \sOX(\vphi^{-1}\Omega_w)_\nu$
  with
  $$
  \sOX(\vphi^{-1}\Omega_w)_\nu = \{f : \vphi^{-1}\Omega_w
  \to \C : f(\bar{g}h)=\nu(h)f(\bar{g}) \; \, \text{for all $h
    \in H$, $\bar{g} \in \vphi^{-1}\Omega_w$}\}.
  $$
  Therefore, by \eqref{L-sheaf}, $\Gamma(\Omega_w,
  \vphi_*\sOX) \cong \bigoplus_{\nu \in \Lambda}
  \Gamma(\Omega_w, \L{\nu})$ as
  $(\g,\cO_\mbB(\Omega_w))$-modules, and it follows that
  $\vphi_*\sOX \cong \bigoplus_{\nu \in \Lambda} \L{\nu}$ as
  $(\g,\cO_\mbB)$-modules.
\end{proof}

One consequence of this proposition is that $\coh^i(\mbX,\sOX)
= 0$ for $i > \dim \mbB = \ell(w_0)=|R^+|$. We can now give the
promised description of $\coh^*(\mbX,\sOX)$.

\begin{thm}
\label{thm57}
As $\DX$ modules, $ \coh^i(\mbX,\sOX) \cong \bigoplus_{w \in
  W(i)}\OXw, $ for $ 0\leq i\leq \dim B$.  Moreover
$$
\coh^*(\mbX,\sOX)\ \cong\ \DX / \DX \g\ \cong\ \bigoplus_{w
  \in W} \OXw
$$
is a direct sum of nonisomorphic simple $\cD$-modules.
\end{thm}

\begin{proof} By Proposition~\ref{cor26} and Lemma~\ref{lem53},
  the $\cO^w$ are nonisomorphic simple modules. Combining
  Theorem~\ref{thm55} with Proposition~\ref{prop56} gives
  $$
  \coh^i(\mbX,\sOX)\ \cong \ \bigoplus_{\nu \in \Lambda}
  \coh^i(\mbB,\L{\nu})
  \  \cong \bigoplus_{\substack{\scriptstyle \mu \in \Lambda^+\\
      \scriptstyle w\in W(i)}} \coh^{i}(\mbB,\L{w\cdot \mu}).
  $$
  Since $\coh^{\ell(w)}(\mbB,\L{w\cdot \mu}) \cong
  V(\mu^*)$, the multiplicity $[\coh^i(\mbX,\sOX) :
  V(\lambda)]$ is equal to $|W(i)|$ for any $\lambda \in
  \Lambda^+$.
  
  Fix $w \in W(i)$ and pick $0\not=e_w$ in the trivial
  $G$-module $\coh^{\ell(w)}(\mbB,\L{w\cdot 0})$.  As $\sOX$ is
  a $(G \times H)$-equivariant $\cD_\mbX$-module, $x e_w =
  \frac{d}{dt}_{\mid t=0}(e^{tx}.e_w) = 0$ for $x \in \g$ and,
  by~\eqref{rem541}, $y \in \wh$ acts on $e_w$ with weight
  $$
  (w\cdot 0)^* = -w_0(w(\rho) - \rho) = w_0ww_0(\rho) -
  \rho.
  $$
  Hence $ye_w = \dual{w_0ww_0(\rho) - \rho}{y} e_w=
  \dual{\rho}{w_0w^{-1}w_0(y) - y} e_w$, for all $y \in \wh$.
  Substituting this into the formula from
  Proposition~\ref{prop22}(2) shows that $ F_{w_0ww_0}(y)
  e_w=0.$
  
  Thus $\big(\cD \g + \cD F_{w_0ww_0}(\wh)\big)e_w = 0$ and so,
  by Theorem~\ref{thm54}(2), $\cD e_w \cong
  \cO^{w_0w^{-1}w_0}$.  Since the map $w \mapsto w_0w^{-1}w_0$
  permutes $W(i)$ and the $\cO^w$ are nonisomorphic simple
  modules, we conclude that $\coh^i(\mbX,\sOX)\supseteq
  \sum_{w\in W(i)} \cD e_w \cong \bigoplus_{w\in W(i)} \cO^w$.
  By Lemma~\ref{g-fixed}, $\cO^w \cong \cO \cong
  \bigoplus_{\lambda \in \Lambda^+} V(\lambda)$ as $G$-modules.
  Thus, for each $\lambda \in \Lambda^+$, the first paragraph
  of the proof implies that
  $$
  \bigl[\ \textstyle{\bigoplus}\{\cO^w : w\in W(i)\} :
  V(\lambda)\bigr]\ =\ |W(i)| \ = \ \bigl[\coh^i(\mbX,\sOX) :
  V(\lambda)\bigr].
  $$
  Consequently, $\coh^i(\mbX,\sOX) \cong \bigoplus_{w \in
    W(i)} \cO^w$.  Theorem~\ref{thm54} then implies that
  $$
  {\coh^*(\mbX,\sOX)\ =\ \bigoplus_{i=0}^{\ell(w_0)}
    \coh^i(\mbX,\sOX)\ \cong\ \bigoplus_{w \in W} \cO^w\ \cong\ 
    \cD / \cD \g,}
  $$
  which completes the proof.
\end{proof}


\section{Differential operators on S-varieties}  
\label{sec3}

In this section we consider highest weight varieties and, more
generally, S-varieties $\ombY$ in the sense of~\cite{VP}. These
are natural generalizations of the closure $\overline{\mbX}$ of
the basic affine space $\mbX$ and there is a natural map from
$\cD(\mbX)$ to the ring of differential operators $\cD(\ombY)$
over such a variety.  Although this map need not be surjective
(see Theorem~\ref{non-exotic}) it does carry enough information
to prove, under mild assumptions, that $\ombY$ is $\cD$-simple
in the sense that $\cO(\ombY)$ is a simple $\cD(\ombY)$-module.
We will continue to write $\cO=\cO(\mbX)$ and $\cD=\cD(\mbX)$.

\begin{defn} \label{Svar} An irreducible affine $G$-variety
  $\overline{\mbY}$ is called an \emph{S-variety} if it
  contains a dense orbit $\mbY= G.v$ such that $U \subseteq
  \mathrm{Stab}_G(v)$, the stabilizer of $v$ in $G$.
\end{defn}

\begin{rem} 
  One important feature of S-varieties is that any affine
  spherical variety (i.e.~an irreducible affine $G$-variety
  having a dense $B$-orbit) is a flat deformation of an
  S-variety (see, for example, \cite[Theorem~22.3]{Gr2}).
\end{rem}

The S-varieties have been completely described in~\cite{VP}.
We begin with the relevant notation.  Set $\Gamma =
\sum_{j=1}^s \N \gamma_j$, where $\gamma_1,\dots,\gamma_s \in
\Lambda^+$ are distinct dominant weights. Write $V_\Gamma =
\bigoplus_{j=1}^s V(\gamma_j)\ni v_\Gamma = v_{\gamma_1} \oplus
\cdots \oplus v_{\gamma_s}$ and define $\mbY_\Gamma =
G.v_\Gamma$.  The following theorem summarizes the results of
\cite[Section~3]{VP} that we need.

\begin{thm}
\label{vinberg-popov}
{\rm (1)} The closures $\overline{\mbY}_\Gamma$ give all the
S-varieties.

\noindent {\rm (2)} One has 
$ \cO(\overline{\mbY}_\Gamma) = \bigoplus_{\gamma \in \Gamma^*}
\cO^{\gamma} \subseteq \cO(\mbY_\Gamma) = \bigoplus_{\mu \in
  \Z\Gamma^* \cap \Lambda^+} \cO^\mu.  $

\noindent {\rm (3)} The following assertions are equivalent:
\begin{enumerate}
\item[(i)] $\cO(\overline{\mbY}_\Gamma)=\cO(\mbY_\Gamma)$;
\item[(ii)] $\Gamma = \Z\Gamma \cap \Lambda^+$;
\item[(iii)] $\overline{\mbY}_\Gamma$ is normal and
  $\codim_{\overline{\mbY}_\Gamma}(\overline{\mbY}_\Gamma
  \smallsetminus \mbY_\Gamma) \ge 2$.\qed
\end{enumerate}
\end{thm}

We will always assume that $\Gamma$ satisfies the equivalent
conditions from Theorem~\ref{vinberg-popov}(3); that is:
\begin{equation}
\label{eq11}
\Z \Gamma \cap \Lambda^+ = \Gamma.
\end{equation}
Hence, for an S-variety $\overline{\mbY}_\Gamma$ satisfying
\eqref{eq11} we have $\cD({\mbY}_\Gamma) =
\cD(\overline{\mbY}_\Gamma)$.  Notice also that $\Gamma^* =
\sum_{j=1}^s \N \gamma_j^*$ satisfies \eqref{eq11} if and only
if $\Gamma$ does.  Natural examples of S-varieties
satisfying~\eqref{eq11} are the following and more examples can
be found in \cite[Ch.~2, \S11]{Gr2}).

\begin{examples}\label{HV-examples}
  {\rm (1)} For $\Gamma = \Lambda^+$ we obtain the basic affine
  space $\mbX = \mbY_{\Lambda^+}$.

\noindent{\rm (2)}  Let $\gamma \in \Lambda^+$ and $\Gamma = \N
\gamma$.  Then, $\mbY_\Gamma$ will be denoted by $\mbY_\gamma$
and is the orbit of a highest weight vector $v_\gamma \in
V(\gamma)$.  Its closure $\overline{\mbY}_\gamma$ is called a
\emph{highest weight} or HV-variety \cite[Section~1]{VP}.
  
\noindent{\rm (3)} An important example
of highest weight variety is the closure of the minimal
(nonzero) orbit in a simple Lie algebra $\g$: in this case
$\gamma = \tilde{\alpha}$ is the highest root and
$V(\tilde{\alpha}) \cong \g$ is the adjoint representation.
\end{examples}

Set $S_\Gamma = \mathrm{Stab}_G(v_\Gamma)$; thus $\mbY_\Gamma
\cong G/S_\Gamma$. Since $S_\Gamma=\bigcap_{j=1}^s
\mathrm{Stab}_G(v_{\gamma_j})$, we see that $\h$ normalizes
$\fs_\Gamma =\mathrm{Lie}(S_\Gamma)$.  Let $\Delta(\fh,
\mathfrak{s}_\Gamma)$ denote the set of roots of $\fh$ in the
Lie algebra $\mathfrak{s}_\Gamma $.  Then $S_\Gamma = S'_\Gamma
Q_\Gamma$ where $Q_\Gamma = H \cap S_\Gamma$ and $S'_\Gamma$ is
generated by the one parameter groups $U_\alpha$, $\alpha \in
\Delta(\fh, \mathfrak{s}_\Gamma)$ (see~\cite[p.~753]{VP}
or~\cite[Corollary~3.5]{Gr2}).  Clearly, $Q_\Gamma = \{ h \in H
: \forall \gamma \in \Gamma, \, \gamma(h) = 1 \}$ is a
diagonalizable subgroup of $H$ with character group
$\Lambda/\Z\Gamma$.  Under the right action of the given groups
on $\cO(\mbY_\Gamma)$, the proof of \cite[Lemma~17.1(b)]{Gr2}
shows that
$$
\cO(\mbY_\Gamma)=\cO(G)^{S_\Gamma} = \cO(G)^{U Q_\Gamma} =
\cO^{Q_\Gamma}.
$$
When \eqref{eq11} holds, this implies that
$\cO(\overline{\mbY}_\Gamma) =\cO(\mbY_\Gamma)= \cO^{Q_\Gamma}$
and $\overline{\mbY}_\Gamma = \overline{\mbX} \twoslash
Q_\Gamma. $

If an algebraic group $L$ acts on the right on some variety
$\mbZ$ we denote by $\delta_g : z \mapsto z.g$ the right
translation by $g \in L$ on $\mbZ$. It induces a right action
on $f \in \cO(\mbZ)$ by $\delta_g.f(z) = f(z.g)$.  This applies
for example to $S_\Gamma$ acting on $G$ and $Q_\Gamma$ acting
on $\overline{\mbX}$.  In this notation, $\cD(\mbX)^\mu = \{d
\in \cD(\mbX) : \forall h \in H, \, \delta_h.d = \mu^*(h) d\}$
for any $\mu \in \Lambda$.

When $\ombY_\Gamma = \ombY_\gamma$ is an HV-variety we will set
$S_\gamma = S_\Gamma$ and $Q_\gamma = Q_\Gamma$, etc.

\begin{prop}
\label{newprop66}
Suppose that $\overline{\mbY}_\Gamma$ is an S-variety
satisfying \eqref{eq11}. Then $\cD^{Q_\Gamma} = \bigoplus_{\mu
  \in \Z\Gamma^*} \cD^{\mu}$ and there is a natural morphism of
algebras $\psi_\Gamma: \cD^{Q_\Gamma} \to \cD(\mbY_\Gamma)$.
\end{prop}

\begin{proof}
  Let $d = \sum_{\mu \in \Lambda} d_{\mu^*}$ with $d_{\mu^*}
  \in \cD(\mbX)^{\mu^*}$. It is clear that $d \in
  \cD(\mbX)^{Q_\Gamma}$ if and only if $\delta_h.d_{\mu^*} =
  d_{\mu^*}$, for all $h \in Q_\Gamma$.  This condition is
  equivalent to $\mu(h)d_{\mu^*}=d_{\mu^*}$; that is, $\mu(h) =
  1$ when $d_{\mu^*} \ne 0$.  Since $Q_\Gamma$ is a
  diagonalizable group, \cite[Proposition~2.5.7(iii)]{Sp}
  implies that $\Z\Gamma = \{\mu \in \Lambda : \forall h \in
  Q_\Gamma, \, \mu(h) = 1\}$. It follows that
  $\cD(\mbX)^{Q_\Gamma} = \bigoplus_{\mu \in \Z\Gamma^*}
  \cD(\mbX)^\mu$. The morphism $\psi_\Gamma$ is simply the
  restriction morphism coming from the identification
  $\overline{\mbY}_\Gamma = \overline{\mbX}\twoslash Q_\Gamma$.
\end{proof}

An affine variety $\mbZ$ (respectively an algebra $R$) is
called {\it $\cD$-simple} if $\cO(\mbZ)$ (respectively $R$) is
a simple $\cD(\mbZ)$-module (respectively $\cD(R)$-module).
This does not hold for arbitrary varieties; for example when
$\mbZ$ is the cubic cone in $\C^3$, $\cO(\mbZ)$ does not even
have finite length as a $\cD(\mbZ)$-module \cite{BGG}.  It is,
however, important in many situations to know that a variety is
$\cD$-simple. We first note that for S-varieties satisfying
\eqref{eq11} this is an easy consequence of
Proposition~\ref{cor26}:

\begin{prop}
\label{thm33'}  
Let $\ombY_\Gamma$ be an S-variety satisfying~\eqref{eq11}.
Then $\cO(\mbY_\Gamma)$ is a simple left
$\cD(\mbY_\Gamma)$-module.
\end{prop}

\begin{proof} 
  Since $\cO(\YGamma)=\OX^{Q_\Gamma}$ with $Q_\Gamma$
  reductive, $\cO(\YGamma)$ is an $\cO(\YGamma)$-module summand
  of $\cO(\mbX)$.  The proposition is now a consequence of
  Proposition~\ref{cor26} and the following result.
\end{proof}

\begin{prop}
\label{prop62}
{\rm (\cite[Proposition~3.1]{Sm})} Let $R \ito T$ be an
inclusion of commutative $\mathbb{C}$-algebras and suppose that
$R$ is a direct summand of the $R$-module $T$. If $T$ is
$\cD$-simple, then $R$ is $\cD$-simple. \qed
\end{prop}

We next refine Proposition~\ref{thm33'} by showing that
$\cO(\mbY_\Gamma)$ is a simple module over a rather explicit
subring of $\cD(\mbY_\Gamma)$.  Let $\calS_\Gamma$ be the
subalgebra of $\DX$ generated by the two finitely generated
commutative subalgebras $\cO(\YGamma)$ and
$F_{w_0}(\cO(\mbY_{\Gamma^*}))$.  Observe that
$\Fo(\cO(\mbY_{\Gamma^*}))$ is isomorphic to
$\cO(\mbY_{\Gamma^*})$ (as both an algebra and a $G$-module)
and so
$$
\calS_\Gamma = \C \langle\, \cO^{\gamma_i^*},
\Fo(\cO^{\gamma_j}) \, ; \, 1 \le i,j \le s \rangle.
$$
By Proposition~\ref{prop22}(1), $\Fo(\cO^\gamma) \subseteq
\cD^{-\gamma^*}$ and so $\calS_\Gamma \subseteq
\cD(\mbX)^{Q_\Gamma}$, by Proposition~\ref{newprop66}. We will
consider $\cO(\YGamma)$ as an $\calS_\Gamma$-module through the
map $\psi_\Gamma$ defined in the latter result.

\begin{prop}
\label{simplicity}
Let $\ombY_\Gamma$ be an S-variety satisfying \eqref{eq11}.
Then $\cO(\YGamma)$ is a simple $\calS_\Gamma$-module.
\end{prop}

\begin{proof}
  As in Lemma~\ref{lem34}, partially order $\Lambda$ by
  $\omega^* \ges \lambda^*$ if $\omega^* - \lambda^* \in
  \Lambda^+$. Let $0 \ne g \in \cO(\YGamma)$ and write $g =
  g_{\lambda_0^*} + \sum_{\lambda_j^* \not\ges \lambda_0^*}
  g_{\lambda_j^*}$ with $g_{\lambda_0^*} \ne 0$ and
  $g_{\lambda_i^*} \in \cO^{\lambda_i^*} \subset \cO(\YGamma)$
  for all $i$. By Lemma~\ref{lem23} there exists $f \in
  \cO^{\lambda_0}$ such that $\Fo(f)(g_{\lambda_0^*})= 1$.
  Hence Proposition~\ref{prop22}(1) implies that $\Fo(f)(g) = 1
  + \sum_{\lambda_j^* \not\ges
    \lambda_0^*}\Fo(f)(g_{\lambda_j^*})$ with
  $\Fo(f)(g_{\lambda_j^*}) \in \cO^{-\lambda_0^* +
    \lambda_j^*}$.  But $\cO^{-\lambda_0^* + \lambda_j^*} =0$
  when $\lambda_j^* \not\ges \lambda_0^*$.  Thus $\Fo(f)(g) =
  1$ and $\cO(\YGamma)$ is simple over $\calS_\Gamma$.
\end{proof}

Let $\mbZ$ be an affine variety with $R=\cO(\mbZ)$ and consider
the subalgebra
$$
\Delta(R)\ =\ \Delta(\mbZ)\ = \ \C\langle\, \cO(\mbZ),
\Der_{\mathbb C}(\cO(\mbZ)\, \rangle\ \subseteq\ \cD(\mbZ).
$$
It is known \cite[Corollary~15.5.6]{MR} that $\Delta(\mbZ) =
\cD(\mbZ)$ when $\mbZ$ is smooth. The Nakai conjecture
\cite{Na} says that the converse should be true:
\begin{equation*}
 \Delta(\mbZ) = \cD(\mbZ) \
 {\buildrel{\text{?}}\over{\limply}} \    R \ \, \text{is
 regular.} 
\end{equation*}
The reader can consult \cite{Be, Tr}, \cite[Section~12.3]{Sc}
and the references therein for work related to this conjecture.

The following observation, which is implicit in \cite{Tr},
implies that many singular varieties, notably S-varieties
satisfying \eqref{eq11}, do satisfy the conclusion of Nakai's
conjecture.

\begin{lem}
\label{prop61}
Assume that $\mbZ$ is irreducible and $\cO(\mbZ)$ is
$\cD$-simple.  Then the Nakai conjecture holds for $\mbZ$.
\end{lem}

\begin{proof}
  Suppose that $\cD(R)=\Delta(R)$. Then $R$ is simple as a
  $\Delta(R)$-module and the result follows from
  \cite[Theorem~15.3.8]{MR}.
\end{proof}

\begin{cor}
\label{cor61}
{\rm (1)} The Nakai conjecture holds for S-varieties satisfying
\eqref{eq11}.

\noindent {\rm (2)} More generally, suppose that $R\subseteq T$
are finitely generated $\mathbb C$-algebras such that $R$ is a
summand of the $R$-module $T$ and that $T$ is $\cD$-simple.
Then the Nakai conjecture holds for $R$.
\end{cor}

\begin{proof} Part (1) is immediate from Lemma~\ref{prop61}
  combined with Proposition~\ref{thm33'}. Similarly, part~(2)
  follows from Lemma~\ref{prop61} and Proposition~\ref{prop62}.
\end{proof}

One significance of the Nakai conjecture is that it implies the
Zariski-Lipman conjecture: $\mbZ$ is smooth whenever $\Der_\C
\cO(\mbZ)$ is a projective $\cO(\mbZ)$-module.  Part~(1) of
Corollary~\ref{cor61} does not give new information about that
conjecture; indeed since S-varieties always have graded
coordinate rings, the Zariski-Lipman conjecture for these
varieties already follows from \cite{Ho}.  It is not clear
whether part~(2) of Corollary~\ref{cor61} has significant
applications in this direction.

A natural situation where Corollary~\ref{cor61} applies is for
invariant rings:

\begin{cor}
\label{thm64}
Let $Q$ be an affine algebraic group and $\mbV$ be an
irreducible affine $Q$-variety. Suppose that $\cO(\mbV)$ is
$\cD$-simple (for example, when $\mbV$ is smooth).  Then the
Nakai conjecture holds for $\mbV \twoslash Q$ in the following
two cases:
\begin{enumerate}
\item[(a)] $Q$ is reductive;
\item[(b)] $\mbV$ is a $G$-variety and $Q=U$.
\end{enumerate}
\end{cor}

\begin{proof}
  (a) As in the proof of Proposition~\ref{thm33'},
  $\cO(\mbV\twoslash Q)$ is a summand of $\cO(\mbV)$. Thus
  Corollary~\ref{cor61}(2) applies.
  
  (b) Observe that $\cD(\overline{\mbX} \times \mbV) \cong
  \cD(\mbX) \otimes \cD(\mbV)$; therefore, by hypothesis and
  Proposition~\ref{cor26} (or Proposition~\ref{simplicity}),
  $\cO(\overline{\mbX} \times \mbV) \cong \OX \otimes
  \cO(\mbV)$ is $\cD$-simple.  By \cite[III.3.2, p.~191]{Kr}),
  $(\overline{\mbX} \times \mbV) \twoslash G \cong \mbV
  \twoslash U$ where $G$ acts componentwise on $\overline{\mbX}
  \times \mbV$.  Thus the result follows from (a) applied to
  $G$ acting on $\overline{\mbX} \times \mbV$.
\end{proof}

Surprisingly, and despite the simplicity of its proof, only
special cases of Corollary~\ref{thm64} have appeared before in
the literature and these have typically required substantially
harder proofs.  See, for example, \cite[Theorem~2.3]{Is} and
\cite[Section~12.3]{Sc}.


\section{Exotic differential operators}
\label{sec4}

The results from the last section raise two questions which we
study in this section. First, can one say more about the
structure, notably the order, of the exotic differential
operators in $\cD(\mbY_\Gamma)$ induced from the ring
$\mathcal{S}_\Gamma$ of Proposition~\ref{simplicity}?  This is
answered by Theorem~\ref{thm44} and Corollary~\ref{exotic-cor}
and proves Corollary~\ref{exotic-intro} from the introduction.
 
The second question concerns the following basic question in
the theory of rings of differential operators. If $V$ is a
finite dimensional representation of a reductive group $K$,
then restriction of operators induces a ring homomorphism
$\cD(V)^K\to \cD(V\twoslash K)$. When is this map surjective?
The conjectural answer is that this happens if and only if
$V\twoslash K$ is singular.  Positive answers to this question
have been found in many cases and these solutions had
significant applications to Lie theory (see, for example,
\cite{Jos2,Jos3, LS0,Sc}). These results have almost always
been in situations where $V\twoslash K$ is a highest weight
variety.  Now let $\overline{\mbY}=\overline{\mbY}_\Gamma$ be
any HV-variety, or even S-variety satisfying \eqref{eq11}.  It
is natural to ask whether the resulting map $\psi_\Gamma
:\cD(\mbX)^{Q_\Gamma} \to \cD(\mbY)$ is surjective.  As we will
show in Theorem~\ref{non-exotic}, this even fails for the
minimal orbit $\mbO_{\mathit{min}}=\mbY_{\tilde{\alpha}}$ of a
simple classical Lie algebra $\g$. This proves
Proposition~\ref{intro-minorbit} from the introduction.

The idea behind the proof of Theorem~\ref{non-exotic} is as
follows.  Let $\ombY=\ombY_\gamma$ be an HV-variety.  Then
$\cO(\mbY) = \bigoplus_{p \in \N} \cO^{p\gamma^*}$ is an
$\N$-graded algebra and so $\cD(\mbY)$ is $\mathbb{Z}$-graded
by \eqref{D-graded}.  There exist examples of HV-varieties
$\ombY$, notably the closures of minimal orbits
$\mbY_{\tilde{\alpha}}$, for which one can find an irreducible
$G$-module $E$ consisting of ``exotic'' operators of degree
$-1$ and order at most $4$, see \cite{AB, BK, LS2}.  On the
other hand, by combining \eqref{D-mu} with
Propositions~\ref{prop22}(1,3) and~\ref{newprop66} the only
obvious operators of order $-1$ in $\Im(\psi_\gamma)$ and small
order are those in $\psi_\gamma(F_{w_0}(\cO^\gamma))$.  But
their order is only bounded above by $k(\gamma)$, where we
write
\begin{equation}\label{k-defn}
\textstyle{k(\lambda) = \dual{\lambda}{2\rho\spcheck} = 2
  \sum_i m_i \qquad\text{for}\quad \lambda =\sum_i m_i\alpha_i
  \in \Lambda^+.} 
\end{equation}
Typically $k(\gamma)$ is significantly larger than $4$ (see the
table at the end of the section).  The aim of the proof is
therefore to show that for the minimal orbit the bound
$k(\gamma)$ is attained and hence that the $G$-module $E$
cannot lie in $\Im(\psi_\gamma)$.
 
We begin with some technical lemmas, the first of which is a
mild generalization of \cite[Lemma~3.8]{BBP}.  The subalgebra
$U(\wh) \cong S(\h)$ of $\cD$ can be identified with
$\cO(\h^*)$ and it follows from the definitions that $u(g)=
u(\lambda)g$ for $u \in U(\wh)$ and $g \in \cO^\lambda$.  We
denote by $\{U_m(\wh)\}_{m \in \N}$ the standard filtration on
the enveloping algebra $U(\wh)$.

\begin{lem}
\label{lem40}
Let $f_1,\dots,f_n \in \cO$ be linearly independent and pick $D
\in \cD_p$.  Suppose that there exist functions $c_i :
\Lambda^+ \to \C$ such that $D(g) = \sum_{i=1}^n c_i(\mu)f_i g$
for all $g \in \cO^\mu$ and $\mu \in \Lambda^+$. Then $D=
\sum_{i=1}^n f_i \tilde{c}_i$ for some $\tilde{c}_1, \dots,
\tilde{c}_n \in U_p(\wh)$.
\end{lem}

\begin{proof}
  For any function $c : \Lambda^+ \to \C$ and $\lambda \in
  \Lambda^+$, define $T_\lambda(c) : \Lambda^+ \to \C$ by
  $T_\lambda(c)(\mu) = c(\lambda + \mu) - c(\mu)$.  Since $\ord
  D \le p$, we have $[g_{p+1},[g_p,[\dots,[g_1,D]] \dots] = 0$
  for all $0 \ne g_j \in \cO^{\lambda_j}$. This easily implies
  that $T_{\lambda_1} \circ T_{\lambda_2} \circ \cdots \circ
  T_{\lambda_p}\circ T_{\lambda_{p+1}}(c_i)=0 $ for all
  $\lambda_i \in \Lambda^+$.  Since $T_\lambda$ is a difference
  operator, it follows that $c_i$ is a polynomial function of
  degree $\le p$ on $\Lambda^+$ and it is clear that there
  exist $\tilde{c}_i \in U(\fh)$ (of degree $\le p$) such that
  $\tilde{c}_i(\lambda) = c_i(\lambda)$ for all $\lambda \in
  \Lambda^+$. Obviously, $D(g) = \left(\sum_{i=1}^n f_i
    \tilde{c}_i\right)(g)$ for all $g \in \cO^\lambda$, hence
  the result.
\end{proof}

\cite[Theorem~1]{Sa1} shows that, for all $\gamma \in
\Lambda^+$ and $ w \in W$, one has
$\cD[\gamma]^{w(\gamma)}\cong U(\wh) \otimes E$, where $E$ is a
$G \times H$-module of dimension $\dim V(\gamma)$.  The next
lemma provides an explicit $G \times H$-module $E$ with this
property.

\begin{lem}
\label{lem41}
Let $\gamma \in \Lambda^+$. Then, via the multiplication map
$\mathsf{m} : \cD \otimes \cD \to \cD$, we have
$$\cD_p[\gamma]^{w(\gamma)} = F_w(\cO^\gamma) \otimes U_p(\wh)
=U_p(\wh) \otimes F_w(\cO^\gamma)$$
for all $w \in W$ and $p
\in \N$.  In particular, $\cD[\gamma]^{w(\gamma)}$ is a free
$U(\wh)$-module with basis being any $\C$-basis of
$F_w(\cO^\gamma)$.
\end{lem}

\begin{proof}
  It is sufficient to prove the lemma for $w=1$; indeed,
  applying $F_w$ to the equalities $\cD_p[\gamma]^{\gamma} =
  \cO^\gamma \otimes U_p(\wh) =U_p(\wh) \otimes \cO^\gamma$ and
  appealing to Proposition~\ref{prop22}(1) gives the general
  result.
  
  By \cite[Theorem~1]{Sa1}, $n=\rk_{U(\wh)} \cD[\gamma]^\gamma=
  \dim V(\gamma)$. Let $\{D_1,\dots, D_n\}$ be a basis of the
  right $U(\wh)$-module $\cD[\gamma]^\gamma$, and fix a basis
  $\{f_1,\dots, f_n\}$ of the $G$-module $\cO^\gamma$.  It is
  easy to see that $U_p(\wh)\otimes \cO^\gamma \cong
  U_p(\wh)\cO^\gamma = \cO^\gamma U_p(\wh) \cong \cO^\gamma
  \otimes U_p(\wh),$ and that this space is contained in
  $\cD[\gamma]^\gamma$.  For the converse, write $f_j = \sum_i
  D_i a_{ij}$ for $a_{ij} \in U(\wh)$. Then, $f_jg = \sum_i
  a_{ij}(\mu) D_i(g)$ for all $0\ne g \in \cO^\mu$.  Thus
  $\bigoplus_{j=1}^n \C f_j g \subseteq \sum_{j=1}^n \C
  D_j(g)$.  Since $\dim \bigl(\bigoplus_{j=1}^n \C f_j g\bigr)
  = n\geq \dim\bigl(\sum_{j=1}^n \C D_j(g) \bigr) $, we obtain
  that $ \bigoplus_{j=1}^n \C f_j g = \bigoplus_{j=1}^n \C
  D_j(g)$.  Therefore, for all $g \in \cO^\mu$, there exist
  unique elements $c_{ij}(g,\mu) \in \C$ such that $D_j(g) =
  \sum_{i=1}^n c_{ij}(g,\mu) f_i g$.
    
  We next show that the $c_{ij}(g,\mu)$'s depend only on $\mu$.
  Indeed, it follows from $$f_j g = \sum_k a_{kj}(\mu) D_k(g)
  =\sum_i\sum_k a_{kj}(\mu) c_{ik}(g,\mu) f_ig$$
  that $\sum_k
  c_{ik}(g,\mu) a_{kj}(\mu)= \delta_{ij}$ for all $1 \le i,j
  \le n$.  In other words, the matrix $[c_{ik}(g,\mu)]_{ik}$ is
  the inverse of the matrix $[a_{kj}(\mu)]_{kj}$. Since the
  $a_{kj}$'s do not depend on $g$, nor do the $c_{ik}(g,\mu)$.
  
  We can therefore write $D_j(g) = {\sum_{i=1}^n} c_{ij}(\mu)
  f_i g$ for any $g \in \cO^\mu$.  By Lemma~\ref{lem40} we
  deduce that $D_j = \sum_{i=1}^n f_i \tilde{c}_{ij}$ with
  $\tilde{c}_{ij} \in U_p(\wh)$.
\end{proof}

Return to an S-variety $\overline{\mbY}=\overline{\mbY}_\Gamma$
satisfying \eqref{eq11} and define $\psi=\psi_\Gamma :
\cD(\mbX)^{Q_\Gamma}\to \cD(\mbY)$ by
Proposition~\ref{newprop66}. Recall that $\Gamma = \sum_{i=1}^s
\N \gamma_i$ for some $\gamma_i\in \Lambda^+$ and order the
$\gamma_i$ so that $\Q \Gamma = \bigoplus_{i=1}^r \Q \gamma_i$.
Set $\ft = \{h \in \wh : \dual{\Gamma^*}{h} = 0\}$; thus $\ft
\cong \Lie(Q_\Gamma)$ has dimension $\text{rank}(\g) -r$.  Pick
$x_1,\dots,x_r \in \wh$ such that $\wh =
\bigl(\bigoplus_{i=1}^r \C x_i\bigr) \oplus \ft$ and
$\dual{\gamma_i^*}{x_j} = \delta_{ij}$ for $1 \le i,j \le r$.
Set $y_j = \psi(x_j) \in \cD(\mbY)$ and let $u \in U(\wh)$.
From the identity $u(f)=\psi(u)(f) = u(\mu)f$ for $\mu \in
\Gamma^*$ and $ f \in \cO^\mu$, one deduces easily that
$\psi(u) = 0$ when $u \in \ft\, U(\wh)$ and that $\psi$ induces
an isomorphism of polynomial algebras:
$$
\psi : \C[\mbx] = \C[x_1,\dots,x_r] \, \longisomto \,
\C[\mby] = \C[y_1,\dots,y_r].
$$
For $\mathbf{m}=(m_1,\dots,m_r) \in \N^r$ we set
$x^\mathbf{m}= \prod_{i=1}^r x_i^{m_i}$ and $y^\mathbf{m} =
\prod_{i=1}^r y_i^{m_i}$.  Note that if $u(y) \in \C[\mby]$,
the (total) degree of $u(y)$ coincides with its order as a
differential operator on $\mbY$.  We recall the definition of
$k(\lambda)$ from \eqref{k-defn}.

\begin{prop}
\label{prop42}  Let $\ombY=\ombY_\Gamma$ 
satisfy \eqref{eq11} and set $\psi=\psi_\Gamma$.  Let $\gamma
\in \Gamma$, $0 \ne f \in \cO^{\gamma}$ and $u(y) \in
\C[\mby]$.  Then,
$$
\ord \psi(\Fo(f)) u(y) = k(\gamma) + \deg u(y).
$$
In particular, the elements $\left\{\psi(\Fo(f_\mathbf{m}))
  y^\mathbf{m} : \mathbf{m} \in \N^r\right\}$ are linearly
independent for any $f_\mathbf{m}\in
\cO^\gamma\smallsetminus\{0\}$.
\end{prop}

\begin{proof}
  As the variety $\overline{\mbY}$ is irreducible, the
  associated graded ring $\gr \cD(\mbY) =\bigoplus_k
  \cD_k(\mbY) / \cD_{k-1}(\mbY)$ is a domain and so $\ord ab
  =\ord a + \ord b$ for $a,b\in \cD(\mbY)$.  Since $\ord u(y) =
  \deg u(y)$, it therefore suffices to show that $\ord
  \psi(\Fo(f)) = k(\gamma)$.  By Proposition~\ref{prop22}(3),
  we do have $\ord \psi(\Fo(f)) \leq \ord \Fo(f)= k(\gamma)$.
  
  Consider the operator $P_\gammas = \sum_i g_i \Fo(f_i)$
  (where $g_i \in \cO^{\gammas}$, $f_i \in \cO^\gamma$) defined
  by Proposition~\ref{prop22}(4).  For each $\alphacheck \in
  \Deltapc$ we have $\alphac = \sum_{j=1}^r
  \dual{\alphac}{\gamma_j^*} x_j + z$ with $z \in \ft$. Thus,
  applying $\psi$ to \eqref{eq-prop22} gives
\begin{equation}
\label{eq44}
\psi(P_\gammas) =   {c_\gammas \prod_{\alphacheck \in
      \Delta_{+}\spcheck}
    \prod_{i=1}^{\dual{\gammas}{\alphacheck}} \bigl(\,\sum_j
    \dual{\alphac}{\gamma_j^*} y_j
    +\dual{\alphacheck}{\rho}-i\bigr).  }
\end{equation}
Write $\gammas = \sum_{j=1}^r m_j \gamma_j^*$ for some $m_j \in
\Q$ and pick $\alpha$ such that $\dual{\gammas}{\alphacheck}
\ne 0$. Then $\sum_j m_j \dual{\alphac}{\gamma_j^*} \ne 0$ and
so $\dual{\alphac}{\gamma_j^*} \ne 0$ for some $1\leq j \leq
r$.  Hence $\deg \psi(\alphac) = \deg\bigl(\sum_j
\dual{\alphac}{\gamma_j^*} y_j \bigr) =1$.  Thus \eqref{eq44}
implies that $\psi(P_\gammas) = \sum_i g_i \psi(\Fo({f_i}))$
has order $\sum_{\alphacheck\in \Delta_{+}\spcheck}
\dual{\gammas}{\alphacheck} = k(\gamma)$. Therefore, there
exists $i \in \{1,\dots,n\}$ such that $\ord \psi(\Fo({f_i})) =
k(\gamma)$. In particular, $\psi(\Fo(\cO^\gamma)) \ne 0$.

Consider the symbol map $\gr_{k(\gamma)} :
\cD_{k(\gamma)}(\mbY) \to \cD_{k(\gamma)}(\mbY) /
\cD_{k(\gamma)-1}(\mbY) $.  Notice that $\gr_{k(\gamma)}$ is a
morphism of $G$-modules and that, by the previous paragraph,
$\gr_{k(\gamma)}(\psi(\Fo({f_i})))\not= 0$ for some $f_i\in
\cO^\gamma$.  Since $\psi(\Fo(\cO^\gamma))\cong V(\gamma)$ is
an irreducible $G$-module, we deduce that
$\gr_{k(\gamma)}(\psi(\Fo({f}))) \ne 0$ for all $0\ne f \in
\cO^\gamma$; that is to say, $\ord \psi(\Fo({f})) = k(\gamma)$.
This proves the first assertion, from which the second follows
immediately.
\end{proof}

\begin{cor}
\label{cor43} Let $\ombY_\Gamma$ satisfy \eqref{eq11}, write 
$\psi=\psi_\Gamma$ and pick $\gamma \in \Gamma$. Then $\psi$
induces an isomorphism of $G$-modules:
$$
\Fo(\cO^\gamma) \otimes \C[\mbx] \, \longisomto \,
\psi(\cD[\gamma]^{-\gammas}) = \psi(\Fo(\cO^\gamma)) \otimes
\C[\mby].
$$
\end{cor}

\begin{proof} We claim that 
  $\cD[\gamma]^{-\gammas}\cap \Ker \psi = \Fo(\cO^\gamma)
  \otimes \ft\, U(\wh)$. Let $Q \in \cD[\gamma]^{-\gammas}\cap
  \Ker \psi$. By Lemma~\ref{lem41}, write $Q= P + T$, where $P
  = \sum_{\mathbf{m}\in \mathbb{N}^r} \Fo(f_\mathbf{m})
  x^\mathbf{m}$ for some $f_\mathbf{m}\in \cO^\gamma$, and $T
  \in \Fo(\cO^\gamma) \otimes \ft\, U(\wh)$. Since $T \in \Ker
  \psi$, certainly $0=\psi(P)= \sum_{\mathbf{m}}
  \psi(\Fo(f_\mathbf{m})) y^\mathbf{m}$.  By
  Proposition~\ref{prop42}, this implies that $f_{\mathbf{m}} =
  0$ for all $\mathbf{m}$.  Hence $P=0$ and $Q \in
  \Fo(\cO^\gamma) \otimes \ft\, U(\wh)$.  Since the opposite
  inclusion is clear, the claim is proven.  As $U(\wh) =
  \C[\mbx] \bigoplus \ft\, U(\wh)$, it follows that $\psi$
  induces the required isomorphism.
\end{proof}

\begin{thm}
\label{thm44}  Let $\ombY_\Gamma$ be an S-variety satisfying
\eqref{eq11}, set $\psi=\psi_\Gamma$ and pick $\gamma \in
\Gamma$. Let $E \subset \psi(\cD[\gamma]^{-\gammas})$ be an
irreducible $G$-module.  Then $E= \psi(\Fo(\cO^\gamma))u(y)$
for some $0\ne u(y) \in \C[\mby]$ and $\ord q = k(\gamma)+\deg
u(y)$ for all $0 \ne q \in E$.

In particular, if $\psi$ is surjective then за$\ord q \ge
k(\gamma)$ for all $0 \ne q \in \cD(\mbY)[\gamma]^{-\gammas}$.
\end{thm}

\begin{proof}
  Let $q \in E \smallsetminus \{0\}$; by Corollary~\ref{cor43}
  we may write $q = \sum_{i} p_i u_i(y)$ where the $p_i \in
  \psi(\Fo(\cO^\gamma))$ are linearly independent and $0 \ne
  u_i(y) \in \C[\mby]$. Since $\psi(\Fo(\cO^\gamma))$ is an
  irreducible $G$-module, the Jacobson Density Theorem produces
  an element $a \in\C G$ such that $ a.p_i = \delta_{i,1}p_1$
  for all $i$.  Recall that $G$ acts trivially on $\wh$, hence
  $a.q = p_1 u_1(y) \in E \smallsetminus \{0\}$. Therefore, $E
  = \C G.(a.q) = \psi(\Fo(\cO^\gamma)) u_1(y)$. This proves the
  first claim and the second is then a consequence of
  Proposition~\ref{prop42}.
  
  Suppose that $\psi$ is surjective. Since $\psi$ is a $(G
  \times H)$-module map, we must have
  $\cD(\mbY)[\gamma]^{-\gammas} =
  \psi(\cD[\gamma]^{-\gamma^*})$ and the final claim follows
  from the previous ones.
\end{proof}

Assume now that $\ombY_\Gamma= \ombY_\gamma$ is an HV-variety
and recall that, by Example~\ref{HV-examples}(2),
$\ombY_\gamma$ automatically satisfies \eqref{eq11}.  Then
$\cO(\mbY_\gamma)$ is graded and from the discussion after
\eqref{D-mu} we know that $\DYg^{-m\gammas}$ identifies with
the space of differential operators of degree $-m$ on the
graded ring $\cO(\Ygamma)$.  In particular, we obtain the
following explicit module of exotic differential operators:

\begin{cor}\label{exotic-cor}
  Let $\ombY=\ombY_\gamma$ be an HV-variety.  Then $\cD(\mbY)$
  contains an irreducible $G$-module
  $E=\psi_\gamma(F_{w_0}(\cO^\gamma))\cong V(\gamma)$ of
  differential operators of degree $-1$ and order $k(\gamma)=
  \dual{\gamma}{2\rho\spcheck}$. \qed
 \end{cor}

 For a number of important HV-varieties, $G$-modules of
 differential operators of degree $-1$ have been constructed,
 but these constructions are typically quite subtle (see, for
 example, \cite{AB,BK0,BK}).  The results of \cite{AB} apply to
 the minimal nilpotent orbit. In our final result we will show
 that their operators almost never appear in the image of
 $\psi$ and therefore that $\psi$ is not surjective for those
 varieties.
 
 Assume that $G$ is simple.  Then the minimal (nonzero)
 nilpotent orbit of $\g=\text{Lie}(G)$ is
 $\mbO_{\mathit{min}}=\mbY_{\tilde{\alpha}}$, where
 $\tilde{\alpha}= \tilde{\alpha}^*$ is the highest root. In
 this case $k(\tilde{\alpha})$ is easy to compute. Indeed, by
 \cite[Ch.~VI, 1.11, Proposition~31]{Bou}, $k(\tilde{\alpha})=
 2(h-1)$, where $h$ is the Coxeter number of the root system
 $\Delta$. These Coxeter numbers are described, for example, in
 \cite[Planche~I--IX]{Bou} and we therefore obtain the
 following values for $k(\tilde{\alpha})$.
 
 \medskip

\begin{center}\begin{tabular}{|c|c|c|c|c|c|c|c|c|c|}
\hline
{\sl Type of $\g$}\vphantom{$\displaystyle\int$} & 
  $\mathsf{A}_\ell$ &
$\mathsf{B}_\ell$, $\ell \ge 2$ &
$\mathsf{C}_\ell$, $\ell \ge 2$ &
$\mathsf{D}_\ell$, $\ell \ge 3$ & 
$\mathsf{E}_6$  &
$\mathsf{E}_7$  &
$\mathsf{E}_8$  & 
$\mathsf{F}_4$  & $\mathsf{G}_2$  \\
\hline {$k(\tilde{\alpha})$} \phantom{$\displaystyle\int$}
 & $ 2\ell$&  $2(2\ell -1)$& $2(2\ell -1)$
& $ 2(2\ell -3)$  & $22$ & $34$& $58$  &  $22$&  $10$\\
\hline
\end{tabular}
\end{center}

\medskip

It is now easy to prove Proposition~\ref{intro-minorbit} from
the introduction.

\begin{thm} \label{non-exotic} 
  Let $\mbO_{\mathit{min}}=\mbY_{\tilde{\alpha}}$ be the
  minimal nonzero orbit in a simple classical Lie algebra $\g$.
  Then the restriction map $\psi :
  \cD(\mbX)^{Q_{\tilde{\alpha}}} \to
  \cD(\mbY_{\tilde{\alpha}})$ is surjective if and only if
  $\g=\fsl(2)$ or $\g=\fsl(3)$.
\end{thm}

\begin{proof} By Example~\ref{HV-examples}(3),
  $\ombY_{\tilde{\alpha}}$ is an HV-variety and so it satisfies
  \eqref{eq11}.

  Suppose that $\psi$ is surjective.  As $\g$ is classical, it
  follows from \cite[Theorem~3.2.3 and Equation~3]{AB} that
  $\cD(\mbY_{\tilde{\alpha}})$ contains a $G$-module $E\cong
  V(\tilde{\alpha})$ of differential operators of degree $-1$
  and order $\leq 4$.  As $\psi$ is a $G\times H$-module map,
  and $\cD(\mbY_{\tilde{\alpha}})^{-\tilde{\alpha}^*}$ is the
  space of operators of degree $-1$, this forces $E\subseteq
  \cD(\mbY_{\tilde{\alpha}})[\tilde{\alpha}]^{-\tilde{\alpha}^*}
  = \psi(\cD[\tilde{\alpha}]^{-\tilde{\alpha}^*})$.  Therefore,
  Theorem~\ref{thm44} says the operators in $E$ have order
  $4\geq k(\tilde{\alpha})$.  The table shows that this is only
  possible when $\g=\fsl(2)$ or $\fsl(3)$.
  
  Conversely, if $\g = \fsl(2)$, then $\cD =
  \C[u,v,\partial_u,\partial_v]$ is the second Weyl algebra,
  $Q_{\tilde{\alpha}} =\{\pm \mathrm{id}\} \cong \Z/2\Z$ and
  $\cO(\mbY_{\tilde{\alpha}}) \ =
  \cO^{Q_{\tilde{\alpha}}}=\C[u^2, v^2, uv]$.  Thus $\psi$ is
  just the isomorphism $\cD^{Q_{\tilde{\alpha}}} = \C[u^2, v^2,
  uv,\partial_u^2,\partial_v^2,\partial_u\partial_v] \isomto
  \cD(\mbY_{\tilde{\alpha}})$.
  
  Now suppose that $\g = \fsl(3)$. Then $\overline{\mbX}$ is
  the quadratic cone $\{\sum_{i=1}^3 u_iy_i = 0\}$ in $\C^3
  \times \C^3$, $Q_{\tilde{\alpha}} \cong \C^*$ and the natural
  map $\psi : \cD^{Q_{\tilde{\alpha}}} \to
  \cD(\mbY_{\tilde{\alpha}})$ is surjective by \cite[Lemma~1.1
  and Theorem~2.14]{LS2}.
\end{proof}

Differential operators have also been extensively studied for
lagran\-gian subvarieties of minimal orbits in \cite{BK0, BK,
  LS2}.  The varieties discussed in those papers are
HV-varieties for an appropriate Lie algebra (see, for example,
\cite[Table~1]{BK0}) and they again have differential operators
of order $\leq 4$ and degree $-1$.  As might be expected by
analogy with Theorem~\ref{non-exotic}, the corresponding map
$\psi$ does produce operators of the required order for Lie
algebras of small rank but in large rank the map $\psi$ is
definitely not surjective. We omit the details of these
assertions since they are rather technical.

The differential operators constructed in \cite{BK0, BK, AB}
have a number of interesting properties, as is explained in
those papers. It would be interesting to know whether the
operators constructed for arbitrary HV-varieties by
Corollary~\ref{exotic-cor} also have distinctive properties.

\begin{rem}\label{lss-comment}
  It is instructive to compare the results of this section with
  those from \cite{LSS,LS0}. One of the main aims of those
  papers was to construct $\cD(\mbZ)$ for certain specific
  singular affine varieties $\mbZ$. The typical situation is
  that $\mbZ$ is an irreducible component of
  $\overline{\mbO\cap \fn^+}$, where $\mbO$ is a nilpotent
  orbit of a simple Lie algebra $\widetilde{\g}$ with
  triangular decomposition $\widetilde{\g}=\fn^-\oplus
  \mathfrak{l} \oplus \fn^+$.  Those papers then show that
  $\cD(\mbZ)=U(\widetilde{\g})/P$ for some primitive ideal $P$.
  However, $\mbZ$ will almost never be an S-variety for
  $\widetilde{\g}$.  Rather, $\mbZ$ will be contained in the
  nilradical of a parabolic subalgebra $\fp$ of
  $\widetilde{\g}$ and, at least when $\mbO$ is the minimal
  orbit, $\mbZ$ then will be an HV-variety for a smaller Lie
  algebra $\g$ contained in the Levi factor of $\fp$.  For
  these examples one would not expect to find a group $Q$ and a
  surjective map $\psi: \cD(G/U)^Q\to \cD(\mbZ)$, simply
  because it is unlikely for the big Lie algebra
  $\widetilde{\g}$ to lie in the image of such a map.
   

  The reader is referred to \cite[Theorem~5.2]{LSS} and
  \cite[Introduction]{LS0} for explicit examples of this
  behaviour and to \cite{Jos2} for a more general framework.
  One example is provided by Example~\ref{sl-3}, for which we
  take $\widetilde{\g}=\mathfrak{so}(8)$. The parabolic $\fp$
  is described explicitly in \cite[Table~3.1 and
  Remark~3.2(v)]{LSS}, so suffice it to say that $\fp$ has
  radical $\fr\cong \C^6$ and Levi factor $\fa\oplus \C$, where
  $\fa\cong \mathfrak{s0}(6) $.  The variety $\mbZ$ is the
  quadric $\sum_{i=1}^3 x_iy_i=0$ inside $\fr$.  Since $\fr$ is
  the natural representation for $\fa$ it follows easily that
  $\mbZ$ is an HV-variety for $\fa$.  Example~\ref{sl-3}
  follows from this discussion by the lucky coincidence that
  the closure $\overline{\mbX}$ of the basic affine space for
  $\mathfrak{sl}(3)$ identifies with $\mbZ$ under an embedding
  of $\mathfrak{sl}(3)$ into $\fa$.
\end{rem}


\end{document}